\newif\ifCPC

\ifCPC 
    \documentclass{cpc}
\else 
    \documentclass[12pt,a4paper,oneside]{article}
    \usepackage{a4}
    \usepackage{amsthm}
    \theoremstyle{plain}
\fi

\usepackage{amsmath, amssymb, bm}
\usepackage{setspace}
\usepackage{xspace}
\usepackage[usenames]{color}
\usepackage{graphicx}
\usepackage[latin1]{inputenc}
\usepackage{multirow}
\usepackage{relsize}
\usepackage[margin=0pt,font=small,labelfont=bf,indention=.45cm]{caption}
\usepackage[indention=0pt,margin=0pt]{subfig}
\usepackage[american]{babel}
\usepackage[numbers,sort&compress,longnamesfirst]{natbib}
\newcommand{\seeurl}[1]{{(also available at {\footnotesize\url{#1}})}}
\newcommand{\toappearseeurl}[1]{{(to appear, see {\footnotesize\url{#1}})}}
\newcommand{\availableurl}[1]{{(available at {\footnotesize\url{#1}})}}

\newcommand{\DIR}{\textsc{dir}}

\newcommand{\INF}{\textsc{inf}}
\newcommand{\CON}{\textsc{con}}
\newcommand{\MINCON}{\textsc{mincon}}
\newcommand{\MAXCON}{\textsc{maxcon}}

\newcommand{\ARR}{\textsc{arr}}
\newcommand{\NEXT}{\textsc{next}}
\newcommand{\EX}{\textsc{ex}}

\newcommand{\V}{V}
\newcommand{\0}{{\ensuremath{ \mathbf{0}}}}

\newcommand{\x}{{\ensuremath{ \mathbf{x}}}}
\newcommand{\y}{{\ensuremath{ \mathbf{y}}}}
\newcommand{\A}{{\ensuremath{ \mathbf{A}}}}
\newcommand{\ee}{\varepsilon}

\newcommand{\fA}{{\ensuremath{ \mathcal{A}}}}
\newcommand{\RS}{{\ensuremath{ \mathcal{R}}}}
\newcommand{\C}{{\ensuremath{ \mathcal{C}}}\xspace}
\newcommand{\B}{{\ensuremath{ \mathcal{B}}}\xspace}

\newcommand{\Aa}{{\ensuremath{ \mathbf{A}^{\!(1)}}}}
\newcommand{\Ab}{{\ensuremath{ \mathbf{A}^{\!(2)}}}}
\newcommand{\Ac}{{\ensuremath{ \mathbf{A}^{\!(3)}}}}
\newcommand{\Ad}{{\ensuremath{ \mathbf{A}^{\!(4)}}}}
\newcommand{\Ax}[1]{{\ensuremath{ \mathbf{A}^{\!(#1)}}}}

\newcommand{\Aba}{{\ensuremath{ \bar{\mathbf{A}}^{\!(1)}}}}
\newcommand{\Abb}{{\ensuremath{ \bar{\mathbf{A}}^{\!(2)}}}}
\newcommand{\Abc}{{\ensuremath{ \bar{\mathbf{A}}^{\!(3)}}}}
\newcommand{\Abd}{{\ensuremath{ \bar{\mathbf{A}}^{\!(4)}}}}
\newcommand{\Aha}{{\ensuremath{ \hat{\mathbf{A}}^{\!(1)}}}}
\newcommand{\Ahb}{{\ensuremath{ \hat{\mathbf{A}}^{\!(2)}}}}
\newcommand{\Ahc}{{\ensuremath{ \hat{\mathbf{A}}^{\!(3)}}}}
\newcommand{\Ahd}{{\ensuremath{ \hat{\mathbf{A}}^{\!(4)}}}}

\newcommand{\Ra}{{\ensuremath{ \mathbf{R}^{\hspace*{-.03\baselineskip}(1)}}}}
\newcommand{\Rb}{{\ensuremath{ \mathbf{R}^{\hspace*{-.03\baselineskip}(2)}}}}
\newcommand{\Rc}{{\ensuremath{ \mathbf{R}^{\hspace*{-.03\baselineskip}(3)}}}}
\newcommand{\Rd}{{\ensuremath{ \mathbf{R}^{\hspace*{-.03\baselineskip}(4)}}}}

\newcommand{\ai}{A^{\hspace*{-.05\baselineskip}(i)} }

\newcommand{\al}{A^{\hspace*{-.05\baselineskip}(\ell)} }
\newcommand{\ax}[1]{A^{\hspace*{-.05\baselineskip}(#1)} }

\newtheorem{thm}{Theorem}
\newtheorem{lem}[thm]{Lemma}
\newtheorem{cor}[thm]{Corollary}

\newcommand{\thmref}[1]{Theorem~\ref{thm:#1}}
\newcommand{\lemref}[1]{Lemma~\ref{lem:#1}}
\newcommand{\lemrefs}[2]{Lemmas~\ref{lem:#1} and~\ref{lem:#2}}

\newcommand{\corref}[1]{Corollary~\ref{cor:#1}}

\newcommand{\assref}[1]{Assumption~\eqref{ass:#1}}

\newcommand{\figref}[1]{Figure~\ref{fig:#1}}

\newcommand{\secref}[1]{Section~\ref{sec:#1}}
\newcommand{\secrefs}[2]{Sections~\ref{sec:#1} and~\ref{sec:#2}}

\newcommand{\eq}[1]{equation~\eqref{eq:#1}}

\newcommand{\Eqs}[2]{Equations~\eqref{eq:#1} and~\eqref{eq:#2}}

\def\argmax{\operatornamewithlimits{argmax}}

\def\mod{\operatorname{mod}}

\DeclareSymbolFont{AMSb}{U}{msb}{m}{n}
\newcommand{\N}{{\mathbb{N}}}

\newcommand{\Z}{{\mathbb{Z}}}
\newcommand{\R}{{\mathbb{R}}}

\newcommand{\arrpp}{\nearrow}
\newcommand{\arrpm}{\searrow}
\newcommand{\arrmp}{\nwarrow}
\newcommand{\arrmm}{\swarrow}
\newcommand{\arru}{\uparrow}
\newcommand{\arrd}{\downarrow}
\newcommand{\arrl}{\leftarrow}
\newcommand{\arrr}{\rightarrow}

\newcommand{\ar}{\ensuremath{\rightarrow}}
\renewcommand{\al}{\ensuremath{\leftarrow}}

\newcommand{\binpp}{\binom{+1}{+1}}
\newcommand{\binpm}{\binom{+1}{-1}}
\newcommand{\binmp}{\binom{-1}{+1}}
\newcommand{\binmm}{\binom{-1}{-1}}

\mathchardef\ordinarycolon\mathcode`\:
\mathcode`\:=\string"8000
\begingroup \catcode`\:=\active
  \gdef:{\mathrel{\mathop\ordinarycolon}}
\endgroup

\allowdisplaybreaks[1]

\ifCPC\else 
    \title{Deterministic Random Walks\\on the Two-Dimensional Grid}
    \author{Benjamin Doerr \and Tobias Friedrich}
    \date{}
\fi

\begin{document}
\ifCPC\label{firstpage}\fi

\maketitle
\begin{abstract}
\noindent
    Jim Propp's rotor router model is a deterministic analogue of a random walk on a 
    graph. Instead of distributing chips randomly, each vertex serves its neighbors 
    in a fixed order. We analyze the difference between Propp machine and random 
    walk on the infinite two-dimensional grid. It is known that, apart from a 
    technicality, independent of the starting configuration, at each time, the 
    number of chips on each vertex in the Propp model deviates from the expected 
    number of chips in the random walk model by at most a constant.
    We show that 
    this constant is approximately 7.8, if all vertices serve their neighbors in 
    clockwise or counterclockwise order and 7.3 otherwise. This result in particular 
    shows that the order in which the neighbors are served makes a difference. Our 
    analysis also reveals a number of further unexpected properties of the two-dimensional
    Propp machine.
\end{abstract}

\ifCPC\else 
    \setstretch{1.1}
\fi

\section{Introduction}

The rotor-router model is a simple deterministic process suggested
by Jim Propp.  It can be viewed as an attempt to derandomize random walks on graphs.
So far, the ``Propp machine'' has mainly been regarded on infinite grids $\Z^d$.
There, each vertex $x\in\Z^d$ is equipped with a ``rotor'' together with a cyclic permutation
(called ``rotor sequence'') of the $2d$ cardinal directions of $\Z^d$.
While a chip (particle, coin, \ldots) performing a random walk leaves a vertex in a random direction,
in the Propp model it always goes into the direction the rotor is pointing.
After a chip is sent, the rotor is rotated according to the fixed rotor sequence.
This shall ensure that the chips are distributed highly evenly among the neighbors.

The Propp machine has attracted considerable attention recently.  It has been shown
that it closely resembles a random walk in several respects.
The first result is due to \citet{Levine1,Levine2} who
compared random walk and Propp machine in an
\emph{aggregating model} called Internal Diffusion-Limited Aggregation (IDLA)~\cite{IDLA}.
There, each chip starts at the origin of $\Z^d$
and walks till it reaches an unoccupied site, which it then occupies.
In the random walk model it is well known that
the shape of the occupied locations converges to a Euclidean ball in $\R^d$~\citep{IDLAlawler1}.
Recently, \citet{Levine1,Levine2} proved an analogous result for the Propp machine.
Surprisingly, the convergence seems to be much faster.
\citet{Kleber} showed experimentally that for circular rotor sequences
after three million chips
the radius of the inscribed and circumscribed circle
differs by approximately $1.61$.
Hence, the occupied locations almost form a perfect circle.
Some more results on this aggregating model in two dimensions can be found in~\secref{blob}.

\citet{CooperComb} compared the Propp machine and the random walk
in terms of the \emph{single vertex discrepancy}. 
Apart from a technicality which we defer to \secref{preliminaries}, they place 
arbitrary numbers of chips on the vertices. Then they run  the Propp machine on 
this initial configuration for a certain number of rounds. A round consists of 
each chip (in arbitrary order) doing one move as directed by the Propp machine. 
For the resulting position, for each vertex they compare the number of chips 
that end up there with the expected number of chips that a random walk in same 
number of rounds would have gotten there starting from the initial 
configuration.
Cooper and Spencer showed that for all grids $\Z^d$, these differences can be 
bounded  by a constant $c_d$ independent of the initial set-up (in particular, 
the total number of vertices) and the run-time.

For the case $d=1$, that is, the graph being the infinite path, \citet{CooperEJC} showed
among other results that this constant $c_1$ is approximately $2.29$.
They further proved that to have the discrepancy on a particular vertex maximal
it suffices that each location has an odd number of chips at at most one time.

In this paper, we rigorously analyze the Propp machine on the two-dimensional grid~$\Z^2$.
A particular difference to the one-dimensional case is that now there are
two non-isomorphic orders in which the four neighbors can be served.
The first are clockwise and counterclockwise orders of the four cardinal directions.
These are called circular rotor sequences.
All other orders turn the rotor by $180^\circ$ at one time and are called
non-circular rotor sequences.
We prove $c_2\approx7.83$ for circular rotor sequences and
$c_2\approx7.29$ otherwise.
To the best of our knowledge, this is the first paper showing that the rotor sequence can make a difference.

We also characterize the respective worst-case configurations.
In particular, we prove that the maximal single vertex discrepancy can 
only be reached if there are vertices which send
a number of chips not divisible by four at at least three different times.

The remainder of this paper is organized as follows. The basic notations
are given in \secref{preliminaries}. In \secref{lower} we roughly speaking show 
that, by suitably choosing the initial configuration, we may prescribe the 
number of chips on each vertex at each time  modulo~4. This will yield sharp 
lower bounds, since in \secref{basic} we see that the discrepancy on a vertex 
can be expressed by exactly this information. 
In \secrefs{INF}{upper}, we derive sufficient information about initial 
configurations leading to maximal discrepancies on a vertex so that we then can 
estimate the maximum possible discrepancy numerically.
This 
estimate is shown to be relatively tight in \secref{error}. Since the 
investigation up to this point in particular showed that different rotor 
sequences lead to different results, we brief{}ly examine the aggregating model 
in this respect in \secref{blob}.
We summarize our results in the last section.


\section{Preliminaries}
\label{sec:preliminaries}
\label{sec:notations}

To bound the single vertex discrepancy between Propp machine and random walk
on the two-dimensional grid we need several preliminaries,
which will be introduced in this section.

First, it will be useful to use a different representation of the two-dimensional grid~$\Z^2$. Let 
$\textstyle\DIR:= \big\{\binpp,\binpm,\binmm,\binmp\big\}$. Define a graph $G=(\V,E)$ via $\V=\big\{\binom{x_1}{x_2} \mid x_1 \equiv x_2\ (\mod 2) \big\}$
and $E=\{ (\x,\y)\in\V^2 \mid \x-\y\in\DIR \}$. Clearly, $G$ is isomorphic to the 
standard two-dimensional grid $G'=(\Z^2,E')$ with
$E'=\{ (\x,\y)\in\Z^2 \mid \|\x-\y\|_1=1 \}$. Therefore, our results on $G$ immediately translate to $G'$.
The advantage of our representation is that now each direction $D \in \DIR$ can 
be uniquely expressed as $D = \ee_x \binom 10 + \ee_y \binom 01$ with $\ee_x, 
\ee_y \in \{-1,1\}$. This allows a convenient computation of the probability 
distribution of the random walk on the grid (see \eq{prob} below).
For convenience we will also use the symbols 
$\big\{\arrpp,\arrpm,\arrmm,\arrmp\big\}$ to describe the directions in the obvious manner.

In order to avoid discussing all equations in the expected sense
and thereby to simplify the presentation, one can treat
the expectation of the random walk as a \emph{linear machine}~\cite{CooperComb}.
Here, in each time step a pile of $k$ chips 
is split evenly, with $k/4$ chips going to each neighbor.
The (possibly non-integral) number of chips
at vertex~$\x$ at time~$t$ is exactly the expected number of chips
in the random walk model.

For $\x,\y\in\V$ and $t\in\N_0$,
let $\x \sim t$ denote that $x_1 \equiv x_2 \equiv t\ (\mod 2)$ and
$\x \sim \y$ denote that $x_1 \equiv x_2 \equiv y_1 \equiv y_2\ (\mod 2)$.
A vertex~$\x$ is called \emph{even} or \emph{odd}
if\, $\x\sim 0$ or $\x\sim 1$, respectively.

A \emph{configuration} describes the current ``state'' of linear or Propp machine.
A configuration of the linear machine assigns to each vertex $\x\in\V$ its current
(possibly fractional) number of chips.
A configuration of the Propp machine assigns to each vertex $\x\in\V$ its current
(integral) number of chips and the current direction of the rotor.
A configuration is called \emph{even} (\emph{odd}) if all chips lie on even
(odd) vertices.

As pointed out in the introduction, there is one limitation without which 
neither the results of \cite{CooperComb,CooperEJC} nor our results hold.  Note that 
since $G$ is a bipartite graph, chips that start on even vertices never 
mix with those starting on odd vertices.  It looks like we are playing 
two games at once.  However, this is not true, because chips at different parity 
vertices may affect each other through the rotors.   We therefore require the initial configuration 
to have chips only on \emph{one} parity. Without loss of generality, we consider 
only even initial configurations.

A random walk on $G$ can be described nicely by its probability density.
By $H(\x,t)$ we denote the probability that a chip from vertex~$\x$
arrives at the origin after~$t$ random steps (``at time~$t$'')
in a simple random walk.
Then,
\begin{equation}
  H(\x,t)=4^{-t} \tbinom{t}{(t+x_1)/2} \tbinom{t}{(t+x_2)/2} \label{eq:prob}
\end{equation}
for $\x\sim t$ and $\|\x\|_\infty\leq t$, and
$H(\x,t)=0$ otherwise.

We now describe the \emph{Propp machine} in detail.
First, we define a rotor sequence by a cyclic permutation $\NEXT\colon \DIR\to\DIR$.
That is, after a chip has been sent in direction~$\A$,
the rotor moves such that afterwards it points in direction $\NEXT(\A)$.
Instead of using $\NEXT$ directly, it will often be more handy to describe
a \emph{rotor sequence} as a
4-tuple $\RS=(\arrpp,\NEXT(\arrpp),\NEXT^2(\arrpp),\NEXT^3(\arrpp))$.
We distinguish between circular and non-circular rotor sequences.
\emph{Circular} rotor sequences
are either clockwise $(\arrpp,\arrpm,\arrmm,\arrmp)$ or
counter-clockwise $(\arrpp,\arrmp,\arrmm,\arrpm)$.
All other rotor sequences are called \emph{non-circular}.
Our main focus is on the classical Propp machine in which all vertices have
the same rotor sequence.
In \citep{CooperComb}, \citeauthor{CooperComb} allow different
rotor sequences for each vertex~$\x$.  Our results also hold in this general setting.
However, to simplify the presentation we will
typically assume that there is only one rotor sequence for all vertices~$\x$.

In the following notations, we implicitly fix the rotor sequence
as well as the initial configuration
(that is, chips on vertices and rotor directions at time $t=0$).
In one step of the Propp machine, each chip does exactly one move, that is, it 
moves in the direction the arrow associated with his current position is pointing 
and updates the arrow direction according to the rotor sequence. Note that the 
particular order in which the chips move within one step is irrelevant (as long 
as we do not label the chips). By this rule, all subsequent configurations are 
determined by the initial configuration. For all $\x \in V$ and $t \in \N_0$ let 
$f(\x,t)$ denote the number of chips on vertex~$\x$ and
$\ARR(\x,t)$ the direction of the rotor associated with~$\x$ after
$t$ steps of the Propp machine.

To describe the linear machine we use the same fixed initial configuration as for the Propp machine.
In one step, each vertex~$\x$ sends a quarter of its (possibly fractional) number of chips
to each neighbor.
Let $E(\x,t)$ denote the number of chips at vertex~$\x$ after~$t$ steps of the linear machine.
This is equal to the expected number of chips at vertex~$\x$
after a random walk of all chips for~$t$ steps.
Note that $E(\x,t)=\tfrac{1}{4} \sum_{\A\in\DIR} E(\x+\A,t-1)$ by definition.


\section{Mod-4-forcing Theorem}
\label{sec:lower}

For a deterministic process like the Propp machine, it is obvious that the 
initial configuration (that is, the location of each chip and the direction of 
each rotor), determines all subsequent configurations. The following theorem 
shows a partial converse, namely that (roughly speaking) we may prescribe the 
number of chips modulo 4 on all vertices at all times and still find an initial 
configuration leading to such a game. An analogous result for the one-dimensional
Propp machine has been shown in \cite{CooperEJC}.

\begin{thm}[Mod-4-forcing Theorem]
    For any initial direction of the rotors
    and any $\pi\colon\V\times\N_0\to\{0,1,2,3\}$ with
    $\pi(\x,t)=0$ for all~$\x\not\sim t$, there is an initial even configuration with
    $f(\x,0)$, $\x\in\V$ that results in subsequent configurations with $f(\x,t)\equiv\pi(\x,t)\ (\mod 4)$
    for all~$\x$ and~$t$.
    \label{thm:parityforcing}
\end{thm}
\begin{proof}
    Let $\ARR(\x,0)$ describe the initial rotor directions given in the assumption. 
    The sought-after configuration can be found iteratively.  We start with
    $f(\x,0):=\pi(\x,0)$ chips 
    at location~$\x$.
    
    Now assume that our initial (even) configuration is such that
    for some $T \in \N$ we have  $f(\x,t)\equiv\pi(\x,t)\ (\mod 4)$ for all~$t<T$.
    We modify this initial configuration by defining
    $f'(\x,0):=f(\x,0)+\varepsilon_\x 4^{T}$ for even~$\x$, while we have
    $f'(\x,0)=0$ for odd~$\x$.  Here,
    $\varepsilon_\x\in\{0,1,2,3\}$ are to be determined such that
    $f'(\x,t)\equiv\pi(\x,t)\ (\mod 4)$ for all~$t\leq T$.
    
    Observe that a pile of $4^{T}$ chips splits evenly~$T$ times. Hence for all 
    choices of the $\varepsilon_\x$ we still have $f'(\x,t)\equiv\pi(\x,t)\ (\mod 4)$
    for all~$t<T$. At time~$T$, the extra piles of $4^T$ chips have spread as 
    follows:
    \[
        f'(\x,T) = f(\x,T) + \!\!\sum_{\substack{\y\sim 0\\ \|\y-\x\|_\infty \leq T}} \!\!\varepsilon_{\y} \binom{T}{\frac{T+x_1-y_1}{2}} \binom{T}{\frac{T+x_2-y_2}{2}}.
    \]
    
    Let initially $\varepsilon_\y:=0$ for all~$\y \in \V$. 
    By induction on $\|\y\|_1$, we change  the $\ee_{\y}$ to their final value.
    We keep $\varepsilon_\y=0$ for all~$\y$ with $\|\y\|_1<2T$.
    
    Assume that for some $\theta \in \N_0$, the current $\varepsilon_\y$ fulfill 
    $f'(\x,T)\equiv\pi(\x,T)\ (\mod 4)$ for all~$\x$ with $\|\x\|_1<\theta$.
    We now determine $\varepsilon_\y$ for all~$\y$ with $\|\y\|_1=2T+\theta$
    in such a way that
    $f'(\x,T)\equiv\pi(\x,T)\ (\mod 4)$ for all~$\x \in \V$ such that $\|\x\|_1\leq\theta$.
    
    Fortunately, to achieve $f'(\x,T)\equiv\pi(\x,T)\ (\mod 4)$ for some $\x \in \V$ 
    such that $\|\x\|_1 = \theta$, it suffices to change a single $\ee_\y$, $\y \in 
    \V$, $\|\y\|_1 = 2T+\theta$. Without loss of generality, let $\x\in\V$, 
    $\|\x\|_1=\theta$, and $\x\sim T$ such that $x_1,x_2\geq 0$. Let $\y = \y(\x) = 
    (x_1+T,x_2+T)$. Now choosing $\varepsilon_{\y} \in \{0,1,2,3\}$ such that 
    $\ee_\y \equiv \pi(\x,T)-f(\x,T)\ (\mod 4)$ yields 
    $f'(\x,T)=f(\x,T)+\varepsilon_{\y} \equiv\pi(\x,T)\ (\mod 4)$ and 
    $f'(\x,T)=f(\x,T)$ for all other~$\x \in \V$ such that $\|\x\|_1 \le \theta$.
    
    Hence for each $\x \in \V$ such that $\|\x\|_1 = \theta$, we find a $\y(\x)$ and 
    a value for $\ee_{\y(\x)}$ such that the resulting $f'(\x,T)$ are as desired. 
    All other $\varepsilon_{\y}$ with $\|\y\|_1=2T+\theta$ remain fixed to zero.
    
    This shows that for all~$\theta \in \N$, we may choose $\varepsilon_\y$, $\y \in 
    \V$, such that $f'(\x,t)\equiv\pi(\x,t)\ (\mod 4)$ for all~$t\leq T$ and all $\x 
    \in \V$ such that $\|\x\|_1\leq\theta$. By compactness principle, this yields 
    the existence of $\varepsilon_\y$, $\y \in \V$, such that 
    $f'(\x,t)\equiv\pi(\x,t)\ (\mod 4)$ for all~$t\leq T$ and $\x \in \V$. 
    
    Up to this point, we proved that for all~$T \in \N$, there is an even initial 
    configuration such that $f'(\x,t)\equiv\pi(\x,t)\ (\mod 4)$ holds for all~$t\leq 
    T$ and $\x \in \V$. Invoking the compactness principle again finishes the proof.
\end{proof}


\section{The Basic Method}
\label{sec:basic}

In this section, we lay the foundations for our analysis of the maximal possible 
single-vertex discrepancy. In particular, we will see that we can determine the 
contribution of a vertex to the discrepancy at another one independent from 
all other vertices.

In the following, we re-use several arguments from \cite{CooperComb,CooperEJC}.
For the moment, in addition to the notations given in \secref{preliminaries},
we also use the following mixed notation.
By $E(\x,t_1,t_2)$ we denote
the (possibly fractional) number of chips at location $x$
after first performing $t_1$ steps with the Propp machine and
then $t_2-t_1$ steps with the linear machine.

We are interested in bounding the discrepancies
$|f(\x,t)-E(\x,t)|$ for all vertices~$\x$ and all times~$t$.
Since we aim at bounds independent of the initial configuration,
it suffices to regard the vertex $\x=\0$.
From
\begin{eqnarray*}
E(\0,0,t) &=& E(\0,t),\\
E(\0,t,t) &=& f(\0,t),
\end{eqnarray*}
we obtain
\begin{equation*}
    f(\0,t)-E(\0,t) \ =\ \sum_{s=0}^{t-1} \left( E(\0,s+1,t) - E(\0,s,t) \right).
\end{equation*}

Now $E(\0,s+1,t) - E(\0,s,t) = \sum_{\x \in V} \sum_{k = 1}^{f(\x,s)}
     \big(H(\x+\NEXT^{k-1}(\ARR(\x,s)), t - s - 1) - H(\x,t-s)\big)$
motivates the definition of the \emph{influence} of a Propp move (compared to a 
random walk move) from vertex~$\x$ in direction $\A$ on the discrepancy of $\0$ 
($t$ time steps later) by
\[
    \INF(\x,\A,t) := H(\x+\A, t-1) - H(\x,t).
\]

To finally reduce all $\ARR$s involved to the initial arrow settings $\ARR(\cdot,0)$, we  define
$s_i(\x) := \min \big\{ u\geq 0 \mid i < \sum_{t=0}^{u} f(\x,t) \big\}$
for all $i\in\N_0$. Hence
at time $s_i(\x)$ the location~$\x$ is 
occupied by its $i$-th chip (where, to be consistent with~\cite{CooperEJC}, we start counting with the $0$-th chip).

Let~$T$ be a time at which we regard the discrepancy at $\0$. Then the above yields
\begin{flalign}\label{eq:basic1}
    f(\0,T) &- E(\0,T) \ =\ \sum_{\x\in \V} \sum_{\substack{i\geq 0,\\ s_i(\x) < T}} \INF(\x,\NEXT^i(\ARR(\x,0)),T-s_i(\x)).
\end{flalign}
Since the inner sum of \eq{basic1} will occur frequently in the remainder, let us define the \emph{contribution} of a vertex~$\x$ to be
    \begin{equation*}
        \CON(\x) \ :=\ \sum_{\substack{i\geq 0,\\ s_i(\x) < T}} \INF(\x,\NEXT^i(\ARR(\x,0)),T-s_i(\x)),
    \end{equation*}
where we both suppress the initial configuration leading to the $s_i(\cdot)$ as well as the run-time~$T$.
Occasionally, we will write $\CON_{\C}$ to specify the underlying initial configuration.

The first main result of this section, summarized in the following theorem,
is that it suffices to examine each vertex~$\x$ separately.
\begin{thm}
    \label{thm:main}
    The discrepancy between Propp machine
    and linear machine after~$T$ time steps is the sum of the contributions $\CON(\x)$
    of all vertices~$\x$,
    i.e.,
    \begin{equation*}
        f(\0,T) - E(\0,T) \ =\ \sum_{\x\in \V} \CON(\x).
    \end{equation*}
\end{thm}
Our  aim in this paper is to prove a sharp upper bound for the single-vertex
discrepancies $|f(\y,T) - E(\y,T)|$ for all~$\y$ and~$T$. As discussed 
already, by symmetry we may always assume $\x = \0$. To get rid of the 
dependency of~$T$, let us define $\MAXCON(\x)$ to be the supremum contribution 
of~$\x$ over all initial configurations and all~$T$. We will shortly see that 
the supremum actually is a maximum (\corref{onlyEX}), that is, there is 
an initial configuration and a time~$T$ such that $\CON(\x) = \MAXCON(\x)$. 
Since the contribution only depends on $T - s_i(\x)$ and the (Mod-4)-forcing theorem
tells us how to manipulate the $s_i(\x)$, we may choose~$T$ as large as 
we like (and still have a configuration leading to $\CON(\x) = \MAXCON(\x)$). 
Provided that $\sum_{\x \in V} \MAXCON(\x)$ is finite (which we prove in the 
remainder), we obtain that $\sum_{\x \in V} \MAXCON(\x)$ is a \emph{tight} upper 
bound for $\sup (f(\0,T) - E(\0,T))$, where the supremum is taken over all 
initial configurations and all~$T$.

To bound $|f(\0,T) - E(\0,T)|$, we need an analogous discussion for negative 
contributions. Let $\MINCON(\x)$ be the infimum contribution of~$\x$ over all 
initial configurations and all $T$. Fortunately, using symmetries, we can show that $\sum_{\x \in V} \MAXCON(\x) = -\sum_{\x \in V} \MINCON(\x)$, hence it suffices to 
regard positive contributions. Let us shortly sketch the symmetry  argument and 
then summarize the above discussion.

Observe that sending one chip in each direction at the same time does not change $\CON(\x)$.
That is, for all~$\x$ and~$t$ we have
\begin{equation}
    \sum_{\A\in\DIR}\!\INF(\x,\A,t) = 0.
    \label{eq:sum0}
\end{equation}
This follows right from the definition of $\INF$ and the elementary fact
$H(\x,t)=\tfrac{1}{4} \sum_{\A\in\DIR} H(\x+\A,t-1)$.
Based on \eq{sum0} 
we will ignore piles of four chips (and multiples) at a common time~$t$
in the remainder of this section.
The remaining one to three chips are called \emph{odd chips}.
Note that there is \emph{no} relation between odd chips and odd vertices/configurations
as defined in \secref{preliminaries}.

To describe the  symmetries of $\CON$, we further distinguish the
non-circular rotor sequences.
\label{alternating}
We call
$(\arrpp,\arrmp,\arrpm,\arrmm)$ and $(\arrpp,\arrmm,\arrpm,\arrmp)$ 
\emph{$x$-alternating} and
$(\arrpp,\arrpm,\arrmp,\arrmm)$ and $(\arrpp,\arrmm,\arrmp,\arrpm)$
\emph{$y$-alternating}.
Now a short look at the definition of $\MAXCON$ reveals symmetries like
$\MAXCON(\binom{x_1}{x_2})=\MAXCON(\binom{-x_1}{-x_2})$
for circular rotor sequences,
$\MAXCON(\binom{x_1}{x_2})=\MAXCON(\binom{x_1}{-x_2})$
for $x$-alternating rotor sequences,
and $\MAXCON(\binom{x_1}{x_2})=\MAXCON(\binom{-x_1}{x_2})$
for $y$-alternating rotor sequences.
The following lemma exhibits symmetries for $\MAXCON$ and $\MINCON$.
It shows that the discrepancies caused by having too few or too many 
chips have the same absolute value. 

\begin{lem}
    \label{lem:sym}
    For all $\x\in V$, the following symmetries hold for
    \begin{itemize}
    \item
    circular rotor sequences:
    $\MAXCON(\binom{x_1}{x_2})=-\MINCON(\binom{-x_1}{x_2})$,
    \item
    $x$-alternating rotor sequences:
    $\MAXCON(\binom{x_1}{x_2})=-\MINCON(\binom{-x_1}{x_2})$,
    \item
    $y$-alternating rotor sequences:
    $\MAXCON(\binom{x_1}{x_2})=-\MINCON(\binom{x_1}{-x_2})$.
    \end{itemize}
\end{lem}
\begin{proof}
    The proofs are not difficult, so we only give the one for the first statement.
    We show that for each  configuration $\C_1$ there is another configuration $\C_3$ and a simple permutation $\pi$ of $V$ with
    $\CON_{\C_1}(\x)= - \CON_{\C_3}(\pi(\x))$ for all implicit run-times $T$ and assuming the clockwise rotor sequence $\RS:=(\arrpp,\arrpm,\arrmm,\arrmp)$ for both $\C_1$ and $\C_3$.
    By \thmref{parityforcing},
    there is a configuration $\C_2$
    which sends,
    using the rotor sequence $(\arrpp,\arrmp,\arrmm,\arrpm)$,
    an odd chip from $\binom{-x_1}{x_2}$ in direction $\binom{-A_1}{A_2}$ at time~$t$
    if and only if $\C_1$ sends an odd chip from~$\binom{x_1}{x_2}$ in direction~$\binom{A_1}{A_2}$ at time~$t$.
    Note that $\CON_{\C_2}(\binom{-x_1}{x_2})=\CON_{\C_1}(\x)$.
    A configuration $\C_3$ which sends for each single chip $\C_2$ sends,
    \emph{three} chips from the same vertex in the same direction at the same time
    obeys rotor sequence $\RS$ and gives
    by \eq{sum0} a contribution
    $\CON_{\C_3}(\binom{-x_1}{x_2}) = 
    -\CON_{\C_2}(\binom{-x_1}{x_2}) = 
    -\CON_{\C_1}(\x)$.
    In consequence, $\MINCON(\binom{-x_1}{x_2})=-\MAXCON(\binom{x_1}{x_2})$ for the clockwise rotor sequence $\RS$.
\end{proof}

Now \lemref{sym} immediately yields $\sum_{\x\in \V} \MINCON(\x) = -\sum_{\x\in \V} \MAXCON(\x)$.
Therefore, it suffices to regard maximal contributions.
\begin{thm}
    \vspace{-\baselineskip}
    \label{thm:tight}
    \[
        \sup_{\C,T} |f(\0,T) - E(\0,T)| \ =\, \sum_{\x\in \V} \MAXCON(\x)
    \]
is a \emph{tight} upper bound for the single
vertex discrepancies.
\end{thm}


\section{The Modes of \smaller\smaller INF}
\label{sec:INF}

In \thmref{tight} we expressed the discrepancy as sum of contributions $\CON(\x)$,
which in turn are sums of the influences $\INF(\x,\A,t)$.
To bound the discrepancy, we are now interested in the extremal values of
such sums.  In this section we derive some monotonicity
properties of these sums.
For this, we define
\[
    \INF(\x,\fA,t) := \sum_{\A\in\fA} \INF(\x,\A,t)
\]
for a finite sequence
$\fA:=(\Aa,\Ab,\ldots)$ of rotor directions ordered according to a
fixed rotor sequence.
In the remainder of the article all finite sequences of
rotor directions for which we use the calligraphic $\fA$
are ordered according to their respective rotor sequence.

\label{def:unimodal}
Let $X \subseteq \R$. We call a
mapping $f\colon X \to \R$ \emph{unimodal}, if there is a $t_1 \in X$ such that
$f|_{x \le t_1}$ as well as $f|_{x \ge t_1}$ are monotone.  We call a
mapping $f\colon X \to \R$ \emph{bimodal}, if there are $t_1,t_2 \in X$ such that
$f|_{x \le t_1}$, $f|_{t_1 \le x \le t_2}$, and $f|_{t_2 \le x}$ are monotone.
We call a mapping $f\colon X \to \R$ \emph{strictly bimodal}, 
if it is bimodal, but not unimodal.
In the following, we show that all $\INF(\x,\fA,t)$ are bimodal in~$t$.

From \eq{sum0} we see that
\begin{equation}
\begin{minipage}[c]{200pt}
\vspace{-.7\baselineskip}
\begin{align}
    \INF(\x,(\Aa,\Ab,\Ac),t) &= 
        -\INF(\x,\DIR\setminus\{\Aa,\Ab,\Ac\},t) \text{\ \ and}\notag\\
    \INF(\x,(\Aa,\ldots,\Ax{k}),t) &= \INF(\x,(\Aa,\ldots,\Ax{k-4}),t) \text{\ \ for $k\geq4$.}\notag
\end{align}
\end{minipage}
    \label{eq:INF0}
\end{equation}
This shows that it suffices to examine $\INF(\x,\fA,t)$ for
$\fA$ of length one and two, which is done in
\lemrefs{INF}{INFF}, respectively.
For both proofs, we need Descartes' Rule of Signs, which can be found
in~\cite{Yap}.
   
\begin{thm}[Descartes' Rule of Signs]
The number of positive roots counting multiplicities
of a non-zero polynomial with real coefficients
is either equal to its number of coefficient sign variations
(i.e., the number of sign changes between consecutive nonzero coefficients)
or else is less than this number by an even integer.
\label{thm:descartes}
\end{thm}

With this, we are now well equipped to analyze the monotonicity of $\INF(\x,\fA,\cdot)$
for $|\fA|\in\{1,2\}$.

\begin{lem}
    For all~$\x\in\V$ and $\A\in\DIR$, $\INF(\x,\A,t)$ is bimodal in~$t$.
    It is strictly bimodal
    if and only if
    \renewcommand{\labelenumi}{(\roman{enumi})}
    \begin{enumerate}
    \item
    $\|\x\|_{\infty}>6$ \ and
    \item
    $-A_1 x_1 > A_2 x_2 > (-A_1 x_1+1)/2$ or
    $-A_2 x_2 > A_1 x_1 > (-A_2 x_2+1)/2$.
    \end{enumerate}
    \label{lem:INF}
\end{lem}
\begin{proof}
    A chip at vertex~$\x$ requires at least $\|\x\|_{\infty}$ time steps to arrive at
    the origin.  Hence, $\INF(\x,\A,t)=0$ for $t< \|\x\|_{\infty}$.
    We show that $\INF(\x,\A,\cdot)$ has at most two
    extrema larger than $\|\x\|_{\infty}$.
    The discrete derivative of $\INF(\x,\A,t)$ in~$t$ is
    \[
        \INF(\x,\A,t+2)-\INF(\x,\A,t) = 
        \frac{p(\x,\A,t) \cdot \big((t-1)!\big)^2 }{4^{t+2}
        \big(\tfrac{t+x_1+2}{2}\big)!
        \big(\tfrac{t-x_1+2}{2}\big)!
        \big(\tfrac{t+x_2+2}{2}\big)!
        \big(\tfrac{t-x_2+2}{2}\big)!
        }
    \]
    with
    $p(\x,\A,t) :=
        (4 A_1 x_1 + 4 A_2 x_2) t^4  +  
        (
        - A_1 x_1^3
        - A_2 x_2^3 
        - A_1 x_1 x_2^2
        - A_2 x_2 x_1^2
        - 6 A_1 x_1 A_2 x_2
        + 19 A_1 x_1
        + 19 A_2 x_2 
        ) t^3  + 
        ( 
          A_1 x_1^3 A_2 x_2 
        + A_1 x_1 A_2 x_2^3 
        - 4 A_1 x_1^3 
        - 4 A_2 x_2^3 
        - 4 A_1 x_1 x_2^2
        - 4 A_2 x_2 x_1^2
        - 23 A_1 x_1 A_2 x_2 
        + 30 A_1 x_1 
        + 30 A_2 x_2
        ) t^2  + 
        (
          A_1 x_1^3 x_2^2
        + A_2 x_2^3 x_1^2
        + 4 A_1 x_1^3 A_2 x_2
        + 4 A_1 x_1 A_2 x_2^3
        - 4 A_1 x_1^3
        - 4 A_2 x_2^3
        - 4 A_1 x_1 x_2^2
        - 4 A_2 x_2 x_1^2
        - 32 A_1 x_1 A_2 x_2
        + 16 A_1 x_1
        + 16 A_2 x_2
        ) t
        - A_1 x_1^3 A_2 x_2^3
        + 4 A_1 x_1^3 A_2 x_2
        + 4 A_1 x_1 A_2 x_2^3
        - 16 A_1 x_1 A_2 x_2.
    $
    We observe that 
    the number of extrema of $\INF(\x,\A,\cdot)$ is exactly the number of roots of 
    $p(\x,\A,\cdot)$.  Since this a polynomial of degree 4 in~$t$, we can use
    Descartes' Sign Rule and some elementary case distinctions to show that 
    $p(\x,\A,\cdot)$ has at most two roots larger than $\|\x\|_{\infty}$.
    A closer calculation reveals that $p(\x,\A,\cdot)$ has precisely two roots larger
    than $\|\x\|_{\infty}$ if
    $\|\x\|_{\infty}>6$ and one of 
    $-A_1 x_1 > A_2 x_2 > (-A_1 x_1+1)/2$ and
    $-A_2 x_2 > A_1 x_1 > (-A_2 x_2+1)/2$ hold.
\end{proof}

\begin{lem}
     For all $\x\in\V$ and $\Aa,\Ab\in\DIR$ such that $\Aa\neq\Ab$, $\INF(\x,(\Aa,\Ab),t)$ is unimodal in~$t$.
    \label{lem:INFF}
\end{lem}
\begin{proof}
    The discrete derivative of $\INF(\x,(\Aa,\Ab),t)$ is
    \begin{eqnarray*}
    \lefteqn{\INF(\x,(\Aa,\Ab),t+2)-\INF(\x,(\Aa,\Ab),t)}\\
    &=& \frac{\big(p(\x,\Aa,t)+p(\x,\Ab,t)\big) \cdot \big((t-1)!\big)^2 }{4^{t+2}
        \big(\tfrac{t+x_1+2}{2}\big)!
        \big(\tfrac{t-x_1+2}{2}\big)!
        \big(\tfrac{t+x_2+2}{2}\big)!
        \big(\tfrac{t-x_2+2}{2}\big)!}
    \end{eqnarray*}
    with $p(\x,\A,t)$ as defined in the proof of \lemref{INF}.
    As there,
    the extrema of $\INF$ are the roots of the quartic function $p(\x,\Aa,t)+p(\x,\Ab,t)$.
    Descartes' Sign Rule now shows that 
    $p(\x,\Aa,t)+p(\x,\Ab,t)$ has at most one root larger than $\|\x\|_{\infty}$ for all~$\x$
    and $\Aa\neq\Ab$.
\end{proof}


\section{Maximal contribution of a vertex}
\label{sec:upper}

We now fix a position~$\x$ and a rotor sequence $\RS$ to examine $\MAXCON(\x)$.
\lemrefs{INF}{INFF} show 
that $\sum_{\A\in\fA} \INF(\x,\A,t)$ is bimodal in~$t$
for all finite sequences
$\fA:=(\Aa,\Ab,\ldots)$ of rotor directions ordered according to $\RS$.
Hence,
for all $\fA$ there are at most two times at which the monotonicity of
$\sum_{\A\in\fA} \INF(\x,\A,t)$ changes.
A time~$t$ at which the monotonicity of 
$\sum_{\A\in\fA} \INF(\x,\A,t)$ changes for some $\fA$
is called \emph{extremal}.
In case of ambiguities, we define the first such time to be extremal.
That is, for unimodal $\sum_{\A\in\fA} \INF(\x,\A,t)$, we choose the first time $t_1$ such that
$\sum_{\A\in\fA} \INF(\x,\A,t)$ is monotone for $t\leq t_1$ and $t\geq t_1$.
Analogously, for strictly bimodal $\sum_{\A\in\fA} \INF(\x,\A,t)$, we choose the first times $t_1$ and $t_2$ such that
$\sum_{\A\in\fA} \INF(\x,\A,t)$ is monotone for $t\leq t_1$, $t_1\leq t\leq t_2$, and $t\geq t_2$.
The set of all \emph{extremal times} is denoted by $\EX(\x)$.

$\EX(\x)$ can be computed easily.
By \eq{INF0} it suffices to consider $\fA$ of length one and two.
The corresponding extremal times are the (rounded) roots of the polynomials 
$p(\x,\A,t)$ and $p(\x,\Aa,t)+p(\x,\Ab,t)$ given in \lemref{INF}.
The following lemma shows that the number of extremal times is very limited.

\begin{lem}
    $|\EX(\x)|\leq 7$
    \label{lem:smallEX}
\end{lem}
\begin{proof}
    According to \lemref{INF}, there is at most one rotor direction~$\A$ for which
    $\INF(\x,\A,t)$ is \emph{strictly} bimodal in~$t$.  Hence, the number of
    extremal times of $\INF(\x,\fA,t)$ with $|\fA|=1$ is at most five.
    For a rotor sequence $\RS=(\Ra,\Rb,\Rc,\Rd)$, \eq{sum0} and
    \lemref{INFF} show that
    $\INF(\x,(\Ra,\Rb),t)=-\INF(\x,(\Rc,\Rd),t)$ and 
    $\INF(\x,(\Rb,\Rc),t)=-\INF(\x,(\Rd,\Ra),t)$ are unimodal in~$t$.
    Therefore, the total number of extremal times of $\INF(\x,(\Aa,\Ab),t)$ with
    $(\Aa,\Ab)$ obeying $\RS$ is at most two.
\end{proof}

Between two successive times $t_1,t_2\in\EX(\x)\cup\{0,T\}$,
$\sum_{\A\in\fA} \INF(\x,\A,t)$ is monotone in~$t$ for all $\fA$.
Such periods of time $[t_1,t_2]$ we call a \emph{phase}.
Note that $\sum_{\A\in\fA} \INF(\x,\A,t)$ could also be constant
in a certain phase.
This implies that it is monotonically increasing as well as
monotonically decreasing.
\label{def:inc}
To avoid this ambiguity, we
use the terms increasing and decreasing (in contrast to
monotonically increasing and decreasing) based
on the minima and maxima at extremal times $\EX(\x)$, which are unambiguously defined and alternating.
We now define precisely when a function 
$\sum_{\A\in\fA} \INF(\x,\A,t)$ is
increasing or decreasing.
Consider the set $E$ of the extremal times of $\sum_{\A\in\fA} \INF(\x,\A,t)$ as defined above.
By \lemrefs{INF}{INFF} we know that $|E|\in\{1,2\}$.
We call $\sum_{\A\in\fA} \INF(\x,\A,t)$ \emph{increasing} at $t$ if
it has a minimum at the maximal $t'\in E$ with $t'< t$ or
a maximum at the minimal $t'\in E$ with $t'> t$.
Analogously, we call $\sum_{\A\in\fA} \INF(\x,\A,t)$ \emph{decreasing} at $t$ if
has a maximum at the maximal $t'\in E$ with $t'< t$ or
a minimum at the minimal $t'\in E$ with $t'> t$.

By abuse of language, let us say that $\x$ sends odd chips at time $t$ if
$f(\x,T-t)\not\equiv 0$ $(\mod 4)$.

\begin{lem}
    \label{lem:onlyEX}
    Let $\C_1$ be an arbitrary configuration with run-time $T\ge \max \EX(\x)$
    and let $\CON_{\C_1}(\x)$ be the corresponding contribution of~$\x$.
    Then there is a configuration $\C_2$ with the same run-time
    and $\CON_{C_1}(\x)\leq\CON_{C_2}(\x)$ that
    sends odd chips only at extremal times, i.e.,
    for the associated $f$ satisfies $f(\x,T-t)\not\equiv 0$ $(\mod 4)$ only
    if $t\in\EX(\x)$.
\end{lem}
\begin{proof}
    Let $\C_2$ be a configuration with
    $\CON_{\C_2}(\x)\geq\CON_{\C_1}(\x)$
    and a minimal number of non-extremal times at which odd chips are sent from~$\x$.
    We assume this number to be greater than zero and show a contradiction.
    
    The sum of the $\INF$s of all chips sent at a certain non-extremal time~$t$
    is either  increasing or  decreasing  in the phase $t$ lies in.
    
    Let us first assume that it is increasing.
    Let $t'$ be the minimal $t'$ such that
    $t'\in\EX(\x)$ or there are odd chips sent at time $t'$
    (assume for the moment that such a $t'$ exists).
    Then, sending the considered pile of odd chips at time $t'$
    instead of time~$t$ decreases the number
    of non-extremal times while not decreasing its contribution.
    Such a modified configuration exists by \thmref{parityforcing}
    and contradicts our assumption on $\C_2$.
    Therefore, there is no such time $t'$.
    This implies that~$t$ lies in the last phase and that the odd chips sent at time~$t$ 
    are the last to be sent at all.
    By $\lim_{t\to\infty}\INF(\x,\A,t)=0$ for all~$\A$, the contribution of
    the chips sent at time~$t$ is negative (since increasing).
    Hence, not sending these chips at all does not decrease $\CON_{\C_2}(\x)$,
    but the number of non-extremal times.
    
    The same line of argument holds
    if the sum of the $\INF$s is decreasing instead of
    increasing.  In this case we use that $\INF(\x,\A,t)=0$ for all
    $t<\|\x\|_\infty$.
\end{proof}

\lemref{onlyEX} immediately gives the following corollary.

\begin{cor}
    There is an initial configuration and a time~$T$ such that $\CON(\x)=\MAXCON(\x)$.
    The configuration can be chosen such that $f(\x,T-t)\not\equiv 0$ $(\mod 4)$ only
    if $t\in\EX(\x)$.
    $T$ can be chosen arbitrarily as long as $T \ge \max \EX(\x)$.
    \label{cor:onlyEX}
\end{cor}

\lemref{smallEX} and \corref{onlyEX} already give a simple, but costly approach to calculate
$\MAXCON(\x)$:  There are four different initial rotor directions for~$\x$ and
at each (of the at most seven) extremal time we can either send 0, 1, 2, or 3 odd chips.
As all subsequent
rotor directions are chosen according to $\RS$, there are only a constant $4\cdot4^7=65536$ configurations to
consider.  The maximum of the respective $\CON(\x)$ will be $\MAXCON(\x)$ by \corref{onlyEX}.

Fortunately, we can also find the worst-case configuration directly.
A \emph{block} of a phase $[t_1,t_2]$ is a 4-tupel $(\Aa,\Ab,\Ac,\Ad)\in\DIR^4$
of rotor directions in the order of $\RS$ such that
$\sum_{i=1}^{k}\INF(\x,\Ax{i},t)$ is increasing in~$t$ in this phase for all~$k\in\{1,2,3\}$.
By \eq{sum0}, this is equivalent to $\sum_{i=k}^{4}\INF(\x,\Ax{i},t)$ being
decreasing in~$t$ within the phase for all~$k\in\{2,3,4\}$.

\begin{lem}
    Each phase has a unique block.  This is determined by the monotonicities of 
    $\INF(\x,\fA,t)$ with $|\fA|\in\{1,2\}$.
    \label{lem:uniq_block}
\end{lem}
\begin{proof}
    Consider a fixed phase.  We want to show that for all valid combinations
    of monotonicities of $\INF(\x,\fA,t)$ with $|\fA|\in\{1,2\}$ within this phase,
    there is exactly one permutation $(\Aa,\Ab,\Ac,\Ad)$ of $\DIR$ obeying $\RS$ such that
    $(\Aa,\Ab,\Ac,\Ad)$ forms a block.
    
    To describe the \emph{type} of monotonicity of $\INF(\x,\A,t)$ within the phase,
    we use a function $\tau$ with $\tau(\A):=\,\ar\,$ if
    $\INF(\x,\A,t)$ is increasing and
    $\tau(\A):=\,\al\,$ if it is decreasing.
    This notation should indicate the direction in which the respective $\INF(\x,\A,t)$ is increasing.
    As a short form we also use $\tau(\Aa,\Ab,\Ac,\Ad):=(\tau(\Aa),\tau(\Ab),\tau(\Ac),\tau(\Ad))$.
    
    By \eq{sum0}, we know that there is at least one~$\A$ of type~$\ar$.
    If there is exactly one direction~$\A$ of type~$\ar$, then the unique permutation
    $(\Aa,\Ab,\Ac,\Ad)$ of $\DIR$ obeying $\RS$ such that
    $\tau(\Aa,\Ab,\Ac,\Ad)=(\ar,\al,\al,\al)$ is the uniquely defined block.
    If there are three rotor directions~$\A$ of type~$\ar$, the block is analogously uniquely defined
    by $\tau(\Aa,\Ab,\Ac,\Ad)=(\ar,\ar,\ar,\al)$.
    
    It remains to examine the case of exactly two rotor directions of type~$\ar$.
    If these two directions are consecutive in $\RS$, 
    $\tau(\Aa,\Ab,\Ac,\Ad)=(\ar,\ar,\al,\al)$ again defines the unique block.
    Otherwise, rotor directions of type~$\ar$ and $\al$ are alternating in the rotor sequence
    and $(\ar,\al,\ar,\al)$ is the only type possible for a block.
    This allows two blocks $(\Aa,\Ab,\Ac,\Ad)$ and $(\Ac,\Ad,\Aa,\Ab)$.
    The choice between these two is uniquely fixed by the monotonicity
    of $\INF(\x,(\Aa,\Ab),t)$.
    Therefore, in all cases there is exactly one unique block.
\end{proof}

We now use \lemref{uniq_block} to define a particular configuration,
which we call \emph{block configuration}.
By \thmref{parityforcing}, we may specify a configuration sufficiently well
by fixing the number of odd chips at all times and locations.
In a block configuration $\B$, a vertex~$\x$ sends odd chips only at extremal times $t\in\EX(\x)$.
Let $(\Aha,\Ahb,\Ahc,\Ahd)$ and $(\Aba,\Abb,\Abc,\Abd)$ denote
the blocks in the phases ending and starting at~$t$.
Then~$\x$ sends $k$~chips at time~$t$ in directions $(\Aa,\ldots,\Ax{k})$,
where $k$ is such that $0\leq k\leq3$ and 
$(\ldots,\Ahd,\Aa,\ldots,\Ax{k},\Aba,\ldots)$ obeys $\RS$.
This uniquely defines when and in which directions odd chips are sent.
Note that we used the blocks only as a technical tool.
There are not necessarily chips sent corresponding to 
$\Aha,\Ahb,\Ahc,\Ahd$ and $\Aba,\Abb,\Abc,\Abd$.
By \thmref{parityforcing}, there are configurations $\B$ as just defined
and for all~$\x$ all of them have the same contribution $\CON_\B(\x)$.


\begin{figure}[bt]
    \centering
    \includegraphics[bb=131pt 157pt 432pt 702pt,angle=-90,width=.75\textwidth,clip]{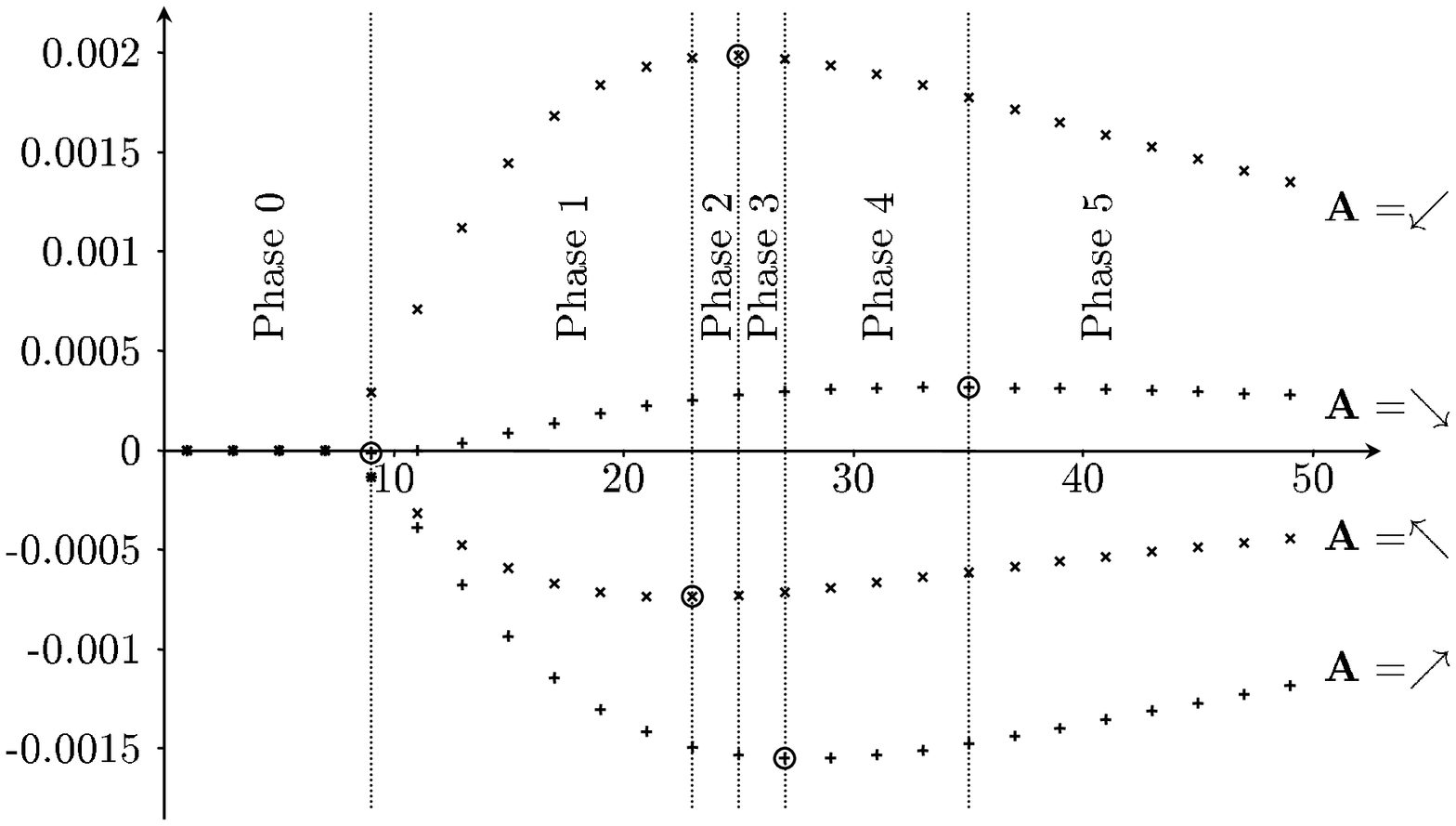}
    \caption{$\INF\big(\binom{5}{9},\A,t\big)$ for $\A\in\{\arrpp,\arrpm,\arrmm,\arrmp\}$.  The circles
             indicate the extrema.}
    \label{fig:INF59}
\end{figure}

\emph{Example.} We now derive the block configuration of the position $\x=\binom{5}{9}$ with the clockwise rotor sequence $\RS=(\arrpp,\arrpm,\arrmm,\arrmp)$.
By calculating the roots of the polynomials $p(\x,\A,t)$ and $p(\x,\Aa,t)+p(\x,\Ab,t)$ given in \lemref{INF},
it is easy to verify that
\begin{itemize}
    \item
    $\INF(\x,\arrpp,t)$ is unimodal with minimum at $t=27$.
    \item
    $\INF(\x,\arrpm,t)$ is bimodal with minimum at $t=9$ and maximum at $t=35$,
    \item
    $\INF(\x,\arrmm,t)$ is unimodal with maximum at $t=25$,
    \item
    $\INF(\x,\arrmp,t)$ is unimodal with minimum at $t=23$,
    \item
    $\INF(\x,(\arrpp,\arrpm),t)$ and
    $\INF(\x,(\arrmp,\arrpp),t)$ are unimodal with minimum at $t=27$.
    \item
    $\INF(\x,(\arrpm,\arrmm),t)$ and
    $\INF(\x,(\arrmm,\arrmp),t)$ are unimodal with maximum at $t=27$.
\end{itemize}
Hence, the extremal points are $\EX(\x)=\{9,23,25,27,35\}$.
\figref{INF59} depicts the plots of $\INF(\x,\A,t)$.
The modes of $\INF(\x,\fA,t)$ listed above uniquely determine the blocks of each phase.
The following table lists rotor directions and type of the block of each phase.
\begin{center}
\begin{tabular}{ccccc}
\hline
\hline
\multirow{2}{*}{\ Phase\ \ } &
\multicolumn{2}{c}{\!\!\!\!\!\!\!\!Boundaries of the phase\ \ \ } &
\multicolumn{2}{c}{\!\!\!\!Block of the phase} \\[-.1cm]
& \ lower & upper & \!\!Rotor directions & Type \\
\hline
0 &  0 &  9 & $\arrmm\arrmp\arrpp\arrpm$ & $\ar\al\al\al$ \\
1 &  9 & 23 & $\arrpm\arrmm\arrmp\arrpp$ & $\ar\ar\al\al$ \\
2 & 23 & 25 & $\arrpm\arrmm\arrmp\arrpp$ & $\ar\ar\al\al$ \\
3 & 25 & 27 & $\arrpm\arrmm\arrmp\arrpp$ & $\ar\al\ar\al$ \\
4 & 27 & 35 & $\arrmp\arrpp\arrpm\arrmm$ & $\ar\ar\ar\al$ \\
5 & 35 & $T$& $\arrmp\arrpp\arrpm\arrmm$ & $\ar\ar\al\al$ \\
\hline\hline
\end{tabular}
\end{center}
This yields the following (maximal as we will see shortly) contribution at $\x=\tbinom{5}{9}$:
\begin{align*}
    \CON(\x)\ =& \ \INF(\x,\arrmm,9)+
                    \INF(\x,\arrmp,9)+
                    \INF(\x,\arrpp,9)+\\
                 & \ \INF(\x,\arrpm,27) + \INF(\x,\arrmm,27) \\
        =& \ \tfrac{20,506,216,364,597}{9,007,199,254,740,992} \approx 0.002277.
\end{align*}
Note that just sending a single chip in the worst direction $\arrmm$
at its worst time $t=25$ gives a smaller contribution of $\INF(\x,\arrmm,25)\approx0.001985$.
Also, sending two chips in directions $\arrpm$ and $\arrmm$ at time
$27=\argmax_t \INF(\x,(\arrpm,\arrmm),t)$ gives
$\INF(\x,(\arrpm,\arrmm),27)\approx0.002261$.
Hence we do profit from sending a chip in the ``wrong'' direction $\arrpp$ at time 9.

The values of $\CON\big(\binom{5}{9}\big)$ for
other rotor sequences are shown in the following table.

    \begin{center}
    \begin{tabular}{ccc}
    \hline
    \hline
    \multirow{2}{*}{Rotor sequence} &
    Times and directions of odd &
    \multirow{2}{*}{$\CON\big(\tbinom{5}{9}\big)$} \\[-.1cm]
    & chips in a block configuration  & \\
    \hline
    $(\arrpp,\arrpm,\arrmm,\arrmp)$ &
    $9:\arrmm\arrmp\arrpp$, $27:\arrpm\arrmm$ &
    \ 0.002277\ldots\!\!
    \\
    $(\arrpp,\arrmp,\arrmm,\arrpm)$ & 
    $23:\arrmm\arrpm\arrpp$, $27:\arrmp\arrmm$, $35:\arrpm$ &
    \ 0.002309\ldots\!\!
    \\
    $(\arrpp,\arrmp,\arrpm,\arrmm)$ &
    $9:\arrmm\arrpp\arrmp$, $23:\arrpm\arrmm\arrpp$, $27:\arrmp\arrpm\arrmm$ &
    \ 0.002302\ldots\!\!
    \\
    $(\arrpp,\arrmm,\arrpm,\arrmp)$ &
    $25:\arrmm$, $35:\arrpm$ &
    \ 0.002230\ldots\!\!
    \\
    $(\arrpp,\arrpm,\arrmp,\arrmm)$ &
    $17:\arrmm\arrpp$, $27:\arrpm\arrmp\arrmm$ &
    \ 0.002083\ldots\!\!
    \\
    $(\arrpp,\arrmm,\arrmp,\arrpm)$ &
    $25:\arrmm$ &
    \ 0.001985\ldots\!\!
    \\
    \hline
    \hline
    \end{tabular}
    \end{center}
    \label{tab:rotors}


\begin{lem}
\label{lem:BC}
    A block configuration yields a contribution of $\MAXCON(\x)$.
\end{lem}
\begin{proof}
    Consider a configuration~$\C$ with contribution $\CON_\C(\x)=\MAXCON(\x)$.
    By previous considerations, we can further assume the following.
    \renewcommand{\labelenumi}{(\arabic{enumi})}
    \begin{enumerate}
    \item
    \label{ass:BC:EX}
    \C only sends odd chips at times $t\in\EX(\x)$
    (cf. \corref{onlyEX}).
    \item
    \label{ass:BC:7}
    \C sends at least seven chips at each time $t\in\EX(\x)$
    (cf. \eq{sum0}).
    \item
    \label{ass:BC:k}
    Let $t_1, t_2 \in \EX(\x)$ such that $[t_1, t_2]$ is a phase and let $k \in \{1,2,3\}$.
    Let $\Ax{1}, \ldots, \Ax{k}$ be the directions the last $k$ chips 
    are sent from vertex~$\x$ at time $t_1$. If 
    $\sum_{i=1}^k\INF(\x,\Ax{i},t)$ is     increasing (cf.~definition on 
    page~\pageref{def:inc}) in $[t_1, t_2]$, then it is not constant. This is a 
    feasible assumption on~$\C$, since otherwise we could send these $k$ 
    chips at time $t_2$ without changing $\CON_\C(\x)$.
    \item
    \label{ass:BC:j}
    Analogously, let $t_1, t_2 \in \EX(\x)$ such that $[t_1, t_2]$ is a phase 
    and let $j \in \{1,2,3\}$. Let $\Ax{1}, \ldots, \Ax{j}$ be the directions the 
    first $j$ chips are sent from vertex~$\x$ at time $t_2$. If 
    $\sum_{i=1}^k\INF(\x,\Ax{i},t)$ is decreasing in $[t_1, t_2]$, then it is not constant. 
    \end{enumerate}
    
    Let $\B$ be a block configuration.
    Aiming at a contradiction, we assume $\CON_\C(\x)>\CON_\B(\x)$.
    Since by \assref{BC:EX} and the definition of \B
    both configurations send odd chips only at times in $\EX(\x)$,
    there is a time $t\in\EX(\x)$ at which the chips of \C contribute
    more than the chips of \B.
    
    We now closely examine the chips sent from~$\x$ at time~$t$ by both 
    configurations. We know that $\B$ sends a uniquely determined number $\ell 
    \in \{0, \ldots, 3\}$ of odd chips at time~$t$ in some directions 
    $\Aa,\ldots,\Ax{\ell}$.
    By the above \assref{BC:7}, $\C$ also sends a sequence of chips in 
    directions $\Aa,\ldots,\Ax{\ell}$. Let $j$ and $k$ denote the number of 
    chips sent by~$\C$ at time~$t$ before and after these $\ell$ chips, 
    respectively. By ignoring, possible piles of four chips, we may assume $j, k 
    \le 3$.
    
    Assume that $k \ge 1$. Then the sum of the $\INF$s of the last $k$ chips
    $\C$ sends at time~$t$ is increasing by the definition of a block. Assume 
    first that~$t$ is not the last extremal time, that is, there is some $t_2 
    \in \EX(\x)$ such that $[t,t_2]$ form a phase. 
    Then by \assref{BC:k} above, the sum of the $\INF$s of the last $k$ chips is 
    strictly increasing in $[t,t_2]$. Hence, a configuration which sends these 
    chips instead at $t_2$ has a larger contribution, in contradiction to the 
    maximality of~$\C$. Now let~$t$ be the last extremal time. From 
    $\lim_{t\to\infty} \INF(\x,\A,t)=0$ for all~$\A$ and the fact that the sum 
    of the $\INF$s of the last $k$ chips is increasing, we see that it is not 
    positive. Hence  the last $k$~chips  do not contribute positively to 
    $\CON_{\C}(\x)$.
        
    Analoguously, assume that $j \ge 1$. Assume first that~$t$ is not the first 
    extremal time, that is, $[t_1,t]$ form a phase for some $t_1 \in \EX(\x)$. 
    By \assref{BC:j}, the first $j$ chips $\C$ sends at time~$t$ have a strictly 
    monotonically decreasing sum of $\INF$s. Hence sending them at time $t_1$ 
    instead of~$t$ gives a larger contribution, again contradicting the 
    maximality of~$\C$. If~$t$ is the first extremal time of~$\x$, then 
    $\INF(\x,\A,t)=0$ for all~$\A$ and $t<\|\x\|_\infty$ shows, similarly as 
    above, that the contribution of the first $j$ chips is not positive.
    
    We conclude that the first $j$ and last $k$ chips sent from~$\x$ and time~$t$
    in~$\C$, if they are present, do not contribute positively to the 
    contribution of~$\x$. This contradicts our assumption 
    $\CON_\C(\x)>\CON_\B(\x)$.
\end{proof}

With the help of a computer, we can now calculate $\MAXCON(\x)$ for all~$\x$.
Using about two months on a Xeon 3 GHz CPU, we computed the maximal contribution
of all vertices in $[-800, 800]^2$. 
If we have the same rotor sequence for all vertices then
\begin{equation}
    \label{eq:res800}
    \sum_{\|\x\|_\infty \leq 800} \MAXCON(\x) = 
                 \begin{cases} 
                    7.832... & \text{for a circular rotor sequence}\\
                    7.286... & \text{for a non-circular rotor sequence.}\\
               \end{cases}\\
\end{equation}
On the other hand, if we allow a different rotor sequence for each vertex, and 
further assume that each vertex has a rotor sequence leading to the maximal 
contribution, then we get
\[
    \sum_{\|\x\|_\infty \leq 800} \MAXCON(\x) = 7.873...
\]
Since all $\MAXCON$s are non-negative,
above values are a lower bound for $\sum_{\x\in V} \MAXCON(\x)$,
and hence for the single vertex discrepancy by \thmref{tight}.

\begin{figure}[p]
    \centering
    \subfloat[Counterclockwise rotor sequence $(\arrpp,\arrmp,\arrmm,\arrpm)$.]{
        \includegraphics[bb=129pt 221pt 480pt 570pt,width=\ifCPC.65\else.6\fi\textwidth,clip]{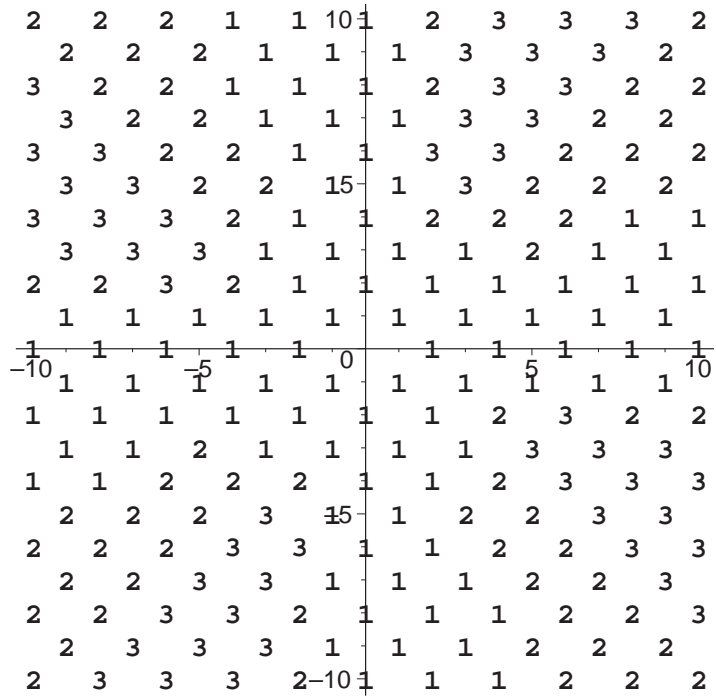}
        \label{subfig:usedpos1}
    }
    \\
    \subfloat[Non-circular rotor sequence $(\arrpp,\arrmp,\arrpm,\arrmm)$.]{
        \includegraphics[bb=129pt 221pt 480pt 570pt,width=\ifCPC.65\else.6\fi\textwidth,clip]{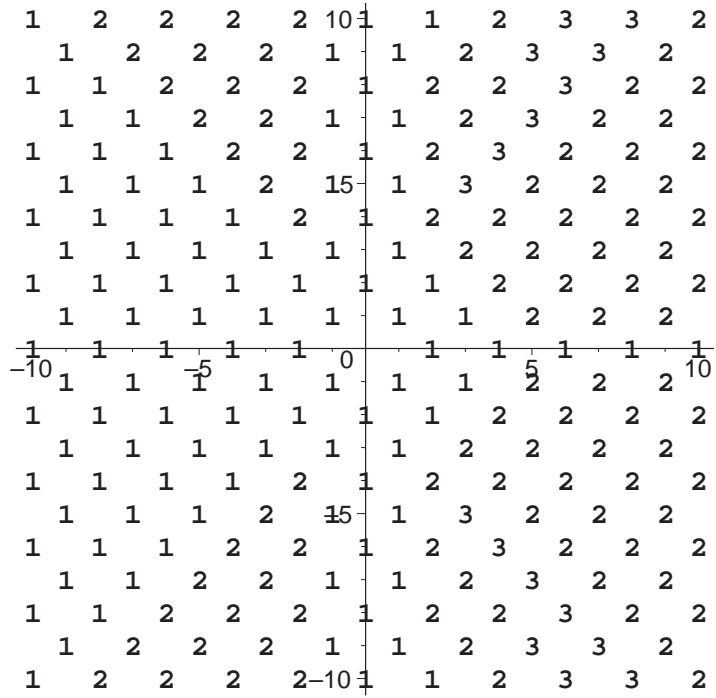}
        \label{subfig:usedpos2}
    }
    \caption{Number of times at which a block configuration sends odd chips.
             Note that the remark on page~\pageref{remark:four} only proved
             $|\widehat\EX(\x)|\leq4$, though there is no position~$\x$
             observable with $|\widehat\EX(\x)|>3$.}
    \label{fig:usedpos}
\end{figure}

\emph{Remark.} 
\label{remark:four}
\lemref{smallEX} shows that the number of extremal times of a vertex is at most seven.
However, a block configuration does not send odd chips at at all extremal times.
Let $\widehat\EX(\x)$ denote the set of extremal times at which odd chips are sent
by the block configuration.
\figref{usedpos} at the end of the paper shows that there are vertices~$\x$ such that $|\widehat\EX(\x)|\geq3$.
We now sketch a proof that $|\widehat\EX(\x)|\leq4$ for all~$\x$.

Note that $|\widehat\EX(\x)|$ only depends on the relative order of the extremal points
and the initial monotonicity (i.e., increasing or decreasing)
of $\INF(\x,\fA,t)$ for $|\fA|\leq 2$.
We use the following two properties of $\INF$ (derived from \eq{sum0}):
\begin{itemize}
\item
In each phase there is at least one $\A\in\DIR$ such that
$\INF(\x,\A,t)$ is increasing (or decreasing).
\item
If $\INF(\x,\Aa,t)$ and $\INF(\x,\Ab,t)$ are both
increasing or decreasing in a phase, so is 
$\INF(\x,(\Aa,\Ab),t)$.
\end{itemize}
For a vertex~$\x$ with only unimodal $\INF(\x,\A,t)$, there are $6!=720$ permutations
of the extrema of $\INF(\x,\A,t)$ and $\INF(\x,(\Aa,\Ab),t)$ and $2^6=64$ initial monotonicities
(using \eq{INF0}).
A simple check by a computer shows that
for only 384 of these 46080 cases both properties from above are satisfied.
For all of them, $|\widehat\EX(\x)|\leq3$ holds.
For vertices~$\x$ with $\INF(\x,\A,t)$ strictly bimodal for an $\A\in\DIR$, there are
$7!/2!=2520$ permutations of the extrema and $2^6=64$ initial monotonicities.
Here, all 408 cases
which satisfy both properties only achieve $|\widehat\EX(\x)|\leq4$.
This proves $|\widehat\EX(\x)|\leq4$ for all~$\x$.
\figref{usedpos} shows $|\widehat\EX(\x)|\leq3$ for all $\|\x\|_\infty \leq 10$.
We could also verify this for $\|\x\|_\infty \leq 800$.
Therefore, we actually expect $|\widehat\EX(\x)|\leq3$ to hold for all~$\x$.
To bridge this gap, stronger properties of $\INF$ seem necessary.


\section{Tail Estimates}
\label{sec:error}

In the previous section, we have calculated the values of
$\sum_{\|\x\|_\infty \leq 800} \MAXCON(\x)$ depending on the rotor sequence.
To show that these are good approximations for the maximal single vertex discrepancy,
we need to find an upper bound on
\[
    E:=\sum_{\|\x\|_\infty > 800} \MAXCON(\x).
\]
In this section, we will prove $E < 0.16$.

We now fix an arbitrary initial configuration and a time $T$.
A simple calculation based on the definitions of $\INF$ and $\CON$ gives
for all~$\x$, $\A$, and $t$
\begin{equation}
\begin{minipage}[c]{200pt}
\vspace{-.7\baselineskip}
\begin{align}
    \INF(\x,\A,t) 
     &= \left(  (A_1 x_1 \cdot A_2 x_2) t^{-2} - (A_1 x_1 + A_2 x_2) t^{-1} \right) H(\x,t),\notag\\
    \CON(\x) &= \sum_{\substack{i\geq 0,\\ s_i(\x) < T}}
                \left(  \frac{\ai_1 x_1 \cdot \ai_2 x_2}{(T-s_i(\x))^{2}} -
                        \frac{\ai_1 x_1 + \ai_2 x_2}{T-s_i(\x)} \right) H(\x,T-s_i(\x))\notag
\end{align}
\end{minipage}
    \label{eq:tail}
\end{equation}
with $s_i(\x)$ as defined in \secref{basic} and $\Ax{i}:=\NEXT^i(\ARR(\x,0))$.
Note that, independent of the chosen rotor sequence,
each of the sequences $(\ai_1(\x))_{i\geq 0}$, $(\ai_2(\x))_{i\geq 0}$, and
$(\ai_1(\x) \ai_2(\x))_{i\geq 0}$ is alternating or alternating in groups of two.
To bound the alternating sums in \eq{tail}, we use
the following fact, which is an elementary extension of Lemma~4 in \cite{CooperEJC}.
\begin{lem}\label{lem:unimodal}
  Let $f\colon X \to \mathbb R$ be non-negative and unimodal with $X\subseteq\R$.  
  Let $\ax{0}, \ldots, \ax{n} \in \{-1,+1\}$ and
  $t_0, \ldots, t_n \in X$ such that $t_0 \leq \ldots \leq t_n$.  
  If $\ai$ is alternating
  or alternating in groups of two, then
    \[\bigg| \sum_{i = 0}^n \ai f(t_i)\bigg| \ \le\  2\max_{x \in X} f(x).\]
\end{lem}

It remains to show that $H(\x,t)/t$ and $H(\x,t) / t^2$ are indeed unimodal.
Note that $\INF(\x,\A,t)$ itself is not always unimodal as shown in \lemref{INF}.

\begin{lem}
    For all~$\x\in V$, $H(\x,t)/t$ and $H(\x,t) / t^2$ are unimodal in~$t$ with
    global maxima at $t_{\max}(\x)$ and $t'_{\max}(\x)$, respectively.
    For the maxima
    $(x_1^2+x_2^2)/4-2 \ \leq\ t_{\max}(\x) \ \leq\ (x_1^2+x_2^2)/4+1$
    and
    $(x_1^2+x_2^2)/6-1 \,\leq\, t'_{\max}(\x) \,\leq\, (x_1^2+x_2^2)/6+2$ holds.
    \label{lem:uni}
\end{lem}
\begin{proof}
    By symmetry, let us assume $x_1 \leq x_2$.
    By definition, $H(\x,t)/t=0$ for $t < x_2$.
    We show that $H(\x,t)/t$ has only one maximum in $t\in[x_2,\infty)$.
    We compute
    \begin{eqnarray*}
    \frac{H(\x,t-2)}{t-2} - \frac{H(\x,t)}{t}
    &=& \frac{4^{-t} p(t) (t-3)!^2\ (t-2)}
         {(\frac{t+x_1}{2})!\ (\frac{t-x_1}{2})!\ (\frac{t+x_2}{2})!\ (\frac{t-x_2}{2})!}.
    \end{eqnarray*}
    with $p(t):=4 t^3 - (x_1^2 + x_2^2 + 5) t^2 + 2 t + x_1^2 x_2^2$.
    By Descartes' Sign Rule (cf. \thmref{descartes}),
    $p(t)$
    has at most one real root larger than $x_2$.
    Since
    \begin{align*}
    p\left(\frac{x_1^2+x_2^2}{4}\right)
        &\ =\ \tfrac{1}{16} \left( 6 x_1^2 x_2^2 + 8 x_1^2 + 8 x_2^2 - 5 x_1^4  - 5 x_2^4 \right)\ <\ 0,\\
    p\left(\frac{x_1^2+x_2^2+5}{4}\right)
        &\ =\ x_1^2 x_2^2 + \frac{x_1^2 + x_2^2 + 5}{2} \ >\ 0.
    \end{align*}
    we see that $H(\x,t)/t$ has a unique extremum, which is a maximum,
    in $[(x_1^2+x_2^2)/4-2,(x_1^2+x_2^2)/4+1]$.
    This proves the lemma for $H(\x,t)/t$.
    The analogous proof for $H(\x,t)/t^2$ is omitted.
\end{proof}
By \eq{tail}, \lemrefs{unimodal}{uni} we obtain
\begin{align*}
E \ &\leq 4 E_1 + 2 E_2
\label{eq:tail2}
\end{align*}
with
\[
    \displaystyle E_1 :=
    \!\sum_{\|\x\|_\infty>800}
            \left| \frac{x_1  H(\x,t_{\max}(\x))}{t_{\max}(\x)}\right|,
    \ \ \ 
    \displaystyle E_2 :=
    \!\sum_{\|\x\|_\infty>800}
            \left| \frac{x_1 x_2 H(\x,t'_{\max}(\x))}{({t'}_{\max}(\x))^2} \right|.
\]
Using \lemref{uni} and $H(\x,t) \leq \big( 2^{-t} \tbinom{t}{t/2} \big)^2 \leq 1/t$,
we now derive upper bounds for $H(\x,t) / t$ and $H(\x,t) / t^2$
for $\|\x\|_{\infty}\ge88$:

\begin{align*}
\left| \frac{H(\x,t_{\max}(\x))}{t_{\max}(\x)} \right|
    \ \leq\ & \frac{1}{t_{\max}(\x)^2} 
    \ \leq\ \frac{16}{(x_1^2+x_2^2-8)^2}
    \ \leq\ \frac{17}{(x_1^2 +x_2^2)^2}, \\
\left| \frac{H(\x,t'_{\max}(\x))}{t'_{\max}(\x)^2} \right|
  \ \leq\ & \frac{1}{t'_{\max}(\x)^3}
  \ \leq\  \frac{216}{(x_1^2+x_2^2-6)^3}  
  \ \leq\  \frac{217}{(x_1^2 +x_2^2)^3}.
\end{align*}
For the calculations in the remainder of this section we need the following estimates.
All of them can be derived by bounding the infinite sums with integrals.
\begin{itemize}
    \addtolength{\itemsep}{-0.1\baselineskip}
    \item
        $\displaystyle\sum_{x>y} \frac{1}{x^k} \leq \frac{1}{(k-1) y^{k-1}}$
        \ \ for all $y > 0$ and all constants $k>1$.
    \item
        $\displaystyle\sum_{x_2=0}^{\infty} \frac{1}{(x_1^2+x_2^2)^2} \leq \frac{7}{3\,x_1^3}$
        \ \ for all $x_1\geq1$.
    \item
        $\displaystyle\sum_{y\geq\beta}
        \frac{1}{(\alpha^2+y^2) y}
        \geq \frac{\ln(\alpha^2+\beta^2)-2 \ln(\beta)}
                {2 \alpha^2}$.
    \item
        $\displaystyle\sum_{\substack{y>\alpha, \\ y \equiv c (\mod 2)}}
            \frac{y}{(y^2+\gamma^2)^2}
            \leq \frac{1}{4(\alpha^2 + \gamma^2)}$.
    \item
        $\displaystyle\sum_{\substack{y>\alpha, \\ y \equiv c (\mod 2)}}
                \frac{1}{(y^2+\gamma^2)^2}
                \leq \frac{\big(\pi - 2 \arctan(\tfrac{\alpha}{\gamma}) \big) (\alpha^2 + \gamma^2)
                          -2 \alpha \gamma}             
                         {8 (\alpha^2+\gamma^2) \gamma^3}$.
    \item        
        $\displaystyle\sum_{y>\beta}
        \frac{\pi - 2 \arctan(\tfrac{\alpha}{y})} {y^2}
        \leq \frac{\ln(\alpha^2+\beta^2) - 2 \ln(\beta)} {\alpha }+
           \frac{\pi - 2 \arctan(\tfrac{\alpha}{\beta})} {\beta}$.
\end{itemize}
With this, we can now bound $E_2$ easily:
\begin{eqnarray}
E_2 &\leq& \!\sum_{\|\x\|_\infty>800}
        \left| \frac{217 x_1 x_2}{(x_1^2 +x_2^2)^3} \right|
    \ \leq\ \!\sum_{\|\x\|_\infty>800}
        \left| \frac{217}{2 (x_1^2 +x_2^2)^2} \right|
        \notag\\        
    &<& \sum_{x_1=1}^{800}
        \sum_{x_2>800}
        \frac{434}{(x_1^2 +x_2^2)^2}
        \ +\
        \!\!\sum_{x_1>800} \sum_{x_2\geq0}
        \frac{434}{(x_1^2 +x_2^2)^2}\notag\\
    &<& \sum_{x_1=1}^{800}
          \sum_{x_2>800}
          \frac{434}{x_2^4}
        \ + \sum_{x_1>800} \frac{3038}{3\,x_1^3}
        \notag\\
    &\leq& \frac{434}{3\cdot 800^2}
        + \frac{1519}{3\cdot 800^2}
    \ < \ 
    0.0011.\label{eq:err2}
\end{eqnarray}
Achieving a good bound for $E_1$ is significantly harder.
We divide $E_1$ in three subsums:
\begin{eqnarray}
    E_1 &<&
        \overbrace{
        4 \sum_{x_1=1}^{800}
        \sum_{\substack{x_2=801, \\ x_2 \equiv x_1 (\mod 2)}}^{\infty}\!\!
        \frac{x_1  H(\x,t_{\max}(\x))}{t_{\max}(\x)}}^\text{see \eq{err1a}}
        \ +
        \ \overbrace{4\!\sum_{x_1=801}^{\infty}
        \sum_{\substack{x_2=1, \\ x_2 \equiv x_1 (\mod 2)}}^{800}\!\!
        \frac{x_1  H(\x,t_{\max}(\x))}{t_{\max}(\x)}}^\text{see \eq{err1b}}\notag\\
        &&
        +\, \underbrace{4\!\sum_{x_1=801}^{\infty}
        \sum_{\substack{x_2=801, \\ x_2 \equiv x_1 (\mod 2)}}^{\infty}\!\!
        \frac{x_1  H(\x,t_{\max}(\x))}{t_{\max}(\x)}}_\text{see \eq{err1c}}
    \ <\ 0.038.\label{eq:err1}
\end{eqnarray}
Now we bound these sums separately as follows.
\begin{eqnarray}
    \lefteqn{4 \sum_{x_1=1}^{800}
               \sum_{\substack{x_2=801, \\ x_2 \equiv x_1 (\mod 2)}}^{\infty}\!\!
               \frac{x_1  H(\x,t_{\max}(\x))}{t_{\max}(\x)}
    \ <\  68 
        \sum_{x_1=1}^{800}
        x_1\!\!\!
        \sum_{\substack{x_2>800, \\ x_2 \equiv x_1 (\mod 2)}}\!\!
        \frac{1}{(x_1^2 +x_2^2)^2} }\notag\\
    &\leq& \frac{17}{2}
        \sum_{x_1=1}^{800}
        \frac{(800^2+x_1^2) \big(\pi-2 \arctan(800/x_1)\big)-1600 \,x_1 }
             {(x_1^2 + 800^2) \, x_1^2} 
    \ <\ 0.0046.\label{eq:err1a}\\
    \lefteqn{4 \sum_{x_1=801}^{\infty}
               \sum_{\substack{x_2=0, \\ x_2 \equiv x_1 (\mod 2)}}^{800}\!\!
               \frac{x_1  H(\x,t_{\max}(\x))}{t_{\max}(\x)}
    \ <\ 68
        \sum_{x_2=0}^{800}  
        \sum_{\substack{x_1>800, \\ x_1 \equiv x_2 (\mod 2)}}\!\!
        \frac{x_1}{(x_1^2 +x_2^2)^2} }\notag\\
    &\leq& 17
        \sum_{x_2=0}^{800}  
        \frac{1}{x_2^2 + 800^2}
    \ <\  0.0167\label{eq:err1b}\\
    \lefteqn{4 \sum_{x_1=801}^{\infty}
               \sum_{\substack{x_2=801, \\ x_2 \equiv x_1 (\mod 2)}}^{\infty}\!\!
               \frac{x_1  H(\x,t_{\max}(\x))}{t_{\max}(\x)}
    \ \leq\ 68
          \sum_{x_1>800}
          x_1\!\!\!
          \sum_{\substack{x_2>800, \\ x_2 \equiv x_1 (\mod 2)}}\!\!
          \frac{1}{(x_1^2 +x_2^2)^2}}\notag\\
    &\leq&\frac{17}{2}
          \sum_{x_1>800}
          \frac{\big(\pi - 2 \arctan(\tfrac{800}{x_1}) \big) (800^2 + x_1^2)-1600 x_1}
                 {(800^2+x_1^2) x_1^2}\notag\\
    &\leq&\frac{17}{2}\left(
          \frac{\ln(2\cdot800^2) - 2 \ln(800)} {800}
          +\frac{\pi/2} {800}
          -\frac{\ln(2\cdot 801^2)-2 \ln(801)} {801}
          \right)\notag\\
    &<&0.0167.\label{eq:err1c}
\end{eqnarray}
Putting this together, we obtain
\begin{align}
E \ &< 4 \cdot 0.038 + 2 \cdot 0.0011 < 0.16.
\label{eq:tail3}
\end{align}
This upper bound on $E$ is not tight. However, it suffices
to prove that the bounds for the single vertex discrepancy calculated
in \secref{upper} do depend on the rotor sequence.
\thmref{tight} and \Eqs{res800}{tail3} yield the following theorem.

\begin{thm}
    \label{thm:res}
    The maximal single vertex discrepancy between Propp machine and linear machine 
    is a constant $c_2$, which depends on the allowed rotor sequences:
    \begin{itemize}
    \item
    If all vertices have the same circular rotor sequence, $7.832\leq c_2\leq 7.985$.
    \item
    If all vertices have the same non-circular rotor sequence, $7.286\leq c_2\leq 7.439$.
    \item
    If all vertices may have different rotor sequences, and we assume
    that each vertex has a rotor sequence leading to a maximal contribution, then $7.873\leq c_2\leq 8.026$.
    \end{itemize}
\end{thm}


\section{Aggregating Model}
\label{sec:blob}

Besides the small single vertex discrepancies examined in the previous sections,
Propp machine and random walk bear striking similarities also in other respects.
The historically first research started by Jim Propp regarded an
aggregating model called Internal Diffusion-Limited
Aggregation (IDLA)~\cite{IDLA}.
In physics this is a well-established model to describe condensation
around a source.

The process starts with an empty grid.
In each round, a particle is inserted at the origin and does a (quasi)-random walk until
it occupies the first empty site it reaches.
For the random walk, it is well known that
the shape of the occupied locations converges to a Euclidean ball~\citep{IDLAlawler1}
in the following sense.
Let $n$ be the number of particles and let
$\Delta(n)$ denote the difference of the radius of the largest inscribed
and the smallest circumscribed circle of an aggregation with $n$ chips.
It has been shown by \citet{IDLAlawler2} that the fluctuations around the limiting shape
are bounded by $\tilde\O(n^{1/6})$ with high probability.
\citet{IDLAmoore} observed experimentally that these error terms
were even smaller, namely poly-logarithmic.

\begin{figure}[p]
    \centering
    \subfloat[Counterclockwise rotor sequence $(\arru,\arrl,\arrd,\arrr)$.]{
        \includegraphics[bb=11pt 10pt 1145pt 1144pt,width=\ifCPC.65\else.6\fi\textwidth,clip]{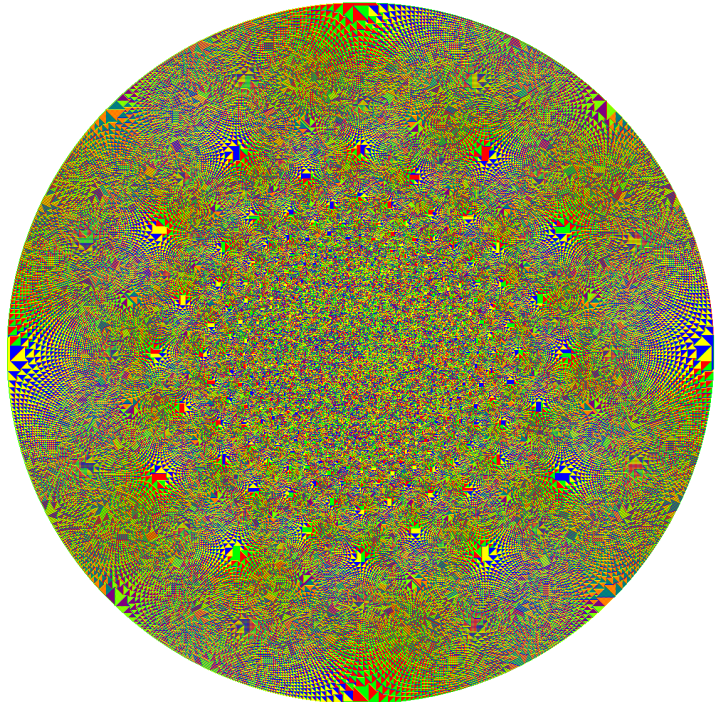}
        \label{subfig:blob1}
    }
    \\
    \subfloat[Non-circular rotor sequence $(\arru,\arrl,\arrr,\arrd)$.]{
        \includegraphics[bb=11pt 10pt 1145pt 1144pt,width=\ifCPC.65\else.6\fi\textwidth,clip]{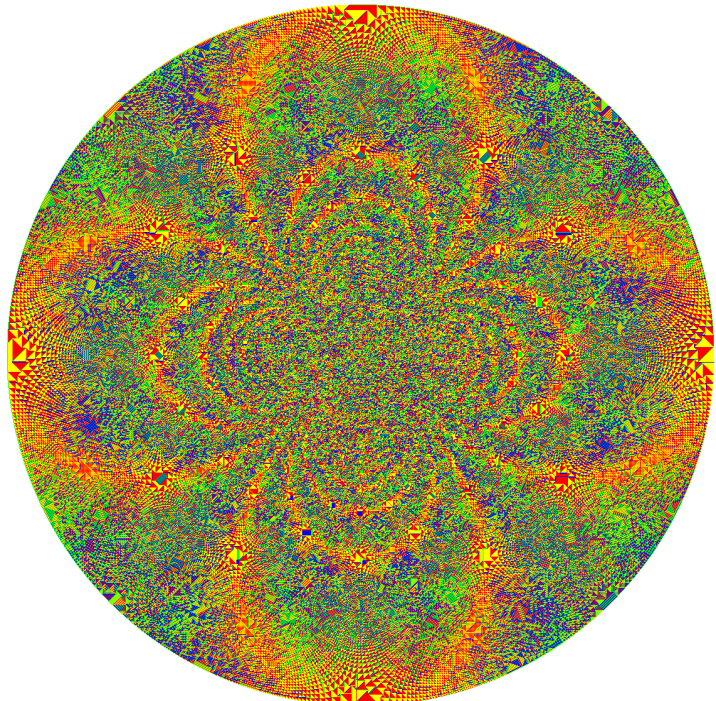}
        \label{subfig:blob2}
    }
    \caption{Propp aggregations with one million particles.
             All rotors initially point to the left.
             The final rotor directions 
             up,    left,   down, and right are denoted by the colors
             red, yellow,  green,  and blue, respectively.}
    \label{fig:blob}
\end{figure}

The analogous model in which the particles do a rotor-router walk instead of a random walk
is much less understood.
\citet{LevineThesis} showed that after $n$ particles have been added,
the Propp aggregation contains a disc of radius $\Omega(n^{1/4})$.
\citet{Levine1,Levine2} proved that 
the shape of occupied locations converges to a Euclidean ball,
however, in a weaker sense than before.  They showed
that the Lebesgue measure of the symmetric difference between the Propp aggregation
and an appropriately scaled Euclidean ball centered at the origin is $O(n^{1/3})$.
Surprisingly, experimental results indicate much stronger bounds.
\citet{Kleber} computed that
for counter-clockwise permutations of the rotor directions 
$\Delta(3\cdot10^6)\approx1.611$
if all rotors initially point to the left.
An apparent conjecture is that there is a constant $\delta$ 
such that $\Delta(n)\leq \delta$ for all $n$.

We reran these experiments with different rotor sequences.
The aggregations for one million particles are shown in \figref{blob}.
Both aggregations do not only differ in the color patterns, but also
in the precise value of $\Delta(n)$.
If all rotors are initially set to the left, we obtained the following values for $\Delta(n)$.
\begin{center}
\begin{tabular}{lccc}
\hline
\hline
Rotor sequence
& $(\arrl,\arru,\arrr,\arrd)$
& $(\arrl,\arru,\arrd,\arrr)$
& $(\arrl,\arrr,\arru,\arrd)$ \\
\hline
average $\Delta(n)$ for $2\cdot10^6<n\leq 3\cdot10^6$ & 1.600  & 0.996 & 1.810 \\
maximal $\Delta(n)$ for $n\leq 3\cdot10^6$ & 1.741 & 1.218 & 1.967 \\
\hline
\hline
\end{tabular}
\end{center}
It is noteworthy that the respective $\Delta(n)$ of both non-circular
rotor sequences $(\arrl,\arru,\arrd,\arrr)$ and
$(\arrl,\arrr,\arru,\arrd)$ differ considerably.

Additionally, we also examined $\Delta(n)$ for random initial rotor directions.
This leads to slightly larger $\Delta$-values.  The following table shows averages and standard deviations
of 100~aggregations with random initial directions of the rotors.
\begin{center}
\begin{tabular}{lccc}
\hline
\hline
Rotor sequence
& $(\arrl,\arru,\arrr,\arrd)$
& $(\arrl,\arru,\arrd,\arrr)$
& $(\arrl,\arrr,\arru,\arrd)$ \\
\hline
average $\Delta(n)$ for $2\cdot10^6<n\leq 3\cdot10^6$ &
    $1.920\pm0.004$ & $1.782\pm0.003$ & $1.781\pm0.003$ \\
maximal $\Delta(n)$ for $n\leq 3\cdot10^6$ &
    $2.541\pm0.051$ & $2.351\pm0.053$ & $2.364\pm0.067$ \\
\hline
\hline
\end{tabular}
\end{center}

As one might have expected, for random initial rotor directions the two non-circular
rotor sequences (columns one and three) are statistically not distinguishable.

The results above again show that different rotor sequences do make a difference.
The main open problem, however, remains to show the conjectured constant upper bound
for~$\Delta(n)$.


\section{Conclusion}
\label{sec:res}

One way of comparing the Propp machine with a random walk is in terms of the 
maximal discrepancy that can occur on a single vertex. It has been shown by 
\citet{CooperComb} that for the underlying graph being an infinite grid~$\Z^d$, 
this single vertex discrepancy can be bounded by a constant $c_d$ independent of the particular initial configuration. For $d=1$, 
this constant has been estimated as $c_1 \approx 2.29$ in~\cite{CooperEJC}. Also, 
the initial configurations leading to a high discrepancy have been described. 
For $d \geq 2$, no such results were known. 

In this paper, we regarded the case $d=2$. We chose the case $d=2$ out of two 
considerations. On the one hand, from dimension two on, there is more than one 
rotor sequence available, which raises the question if different rotors 
sequences make a difference. One the other hand, we restrict ourselves to $d=2$, 
because for larger $d$ a nice expression for the probability $H(\x,t)$ that a 
chip from vertex~$\x$ arrives at the origin after $t$ random steps is missing. 
This probably makes it very hard to find sufficiently sharp estimates for the 
single vertex discrepancies.

We were able to give relatively tight estimates for the constants $c_2$ taking 
into account different rotor sequences and obtain several interesting facts 
about the worst-case initial configurations.
The maximal single vertex discrepancy $c_2$ satisfies the following.
If all vertices have the same circular rotor sequence, $7.832\leq c_2\leq 7.985$.
If all vertices have the same non-circular rotor sequence, $7.286\leq c_2\leq 7.439$.
If all vertices may have different rotor sequences, and we assume
that each vertex has a rotor sequence leading to a maximal contribution, then $7.873\leq c_2\leq 8.026$.
In particular, we see that non-circular rotor sequences seem to produce smaller 
discrepancies than circular one. The gaps between upper and lower bounds stem 
from the fact that we used a computer to calculate the precise maximal 
contribution $\CON(\x - \y)$ of vertex~$\x$ on the discrepancy at $\y$. Hence 
the lower bounds are the maximal discrepancies obtained from initial 
configurations such that all vertices~$\x$ with $\|\x-\y\|_{\infty} > 800$  at all times 
contain numbers of chips only that are divisible by 4.

We also learned that the initial configurations leading to such 
discrepancies are more complicated than in the one-dimensional
case. Recall from~\cite{CooperEJC} that in the one-dimensional case in 
a worst-case setting each position needs to have an odd number of chips only 
once. If we aim at a surplus of chips in the Propp model, these odd chips were 
always sent towards the position under consideration, otherwise away from it.

In the two-dimensional case, things are more complicated. Here it can be 
necessary that a position holds a number of chips not divisible by 4 up to three 
times.
Also, the number of ``odd'' chips (those which cannot be put into piles of four) 
can be as high as nine. In consequence, it can make sense to send odd chips in 
the wrong direction (e.g., away from the position where we aim at a surplus of 
chips). An example showing this was analyzed in Section~\ref{sec:upper}. The 
reason for such behavior seems to be that the influences $\INF(\x,A,t)$ of odd 
chips sent from~$\x$ in direction $A$ at time $t$ are not unimodal functions in 
$t$ anymore (as in the one-dimensional case).

\begin{figure}[p]
    \centering
    \subfloat[Counterclockwise rotor sequence $(\arrpp,\arrmp,\arrmm,\arrpm)$.]{
        \includegraphics[bb=129pt 221pt 480pt 570pt,width=\ifCPC.65\else.6\fi\textwidth,clip]{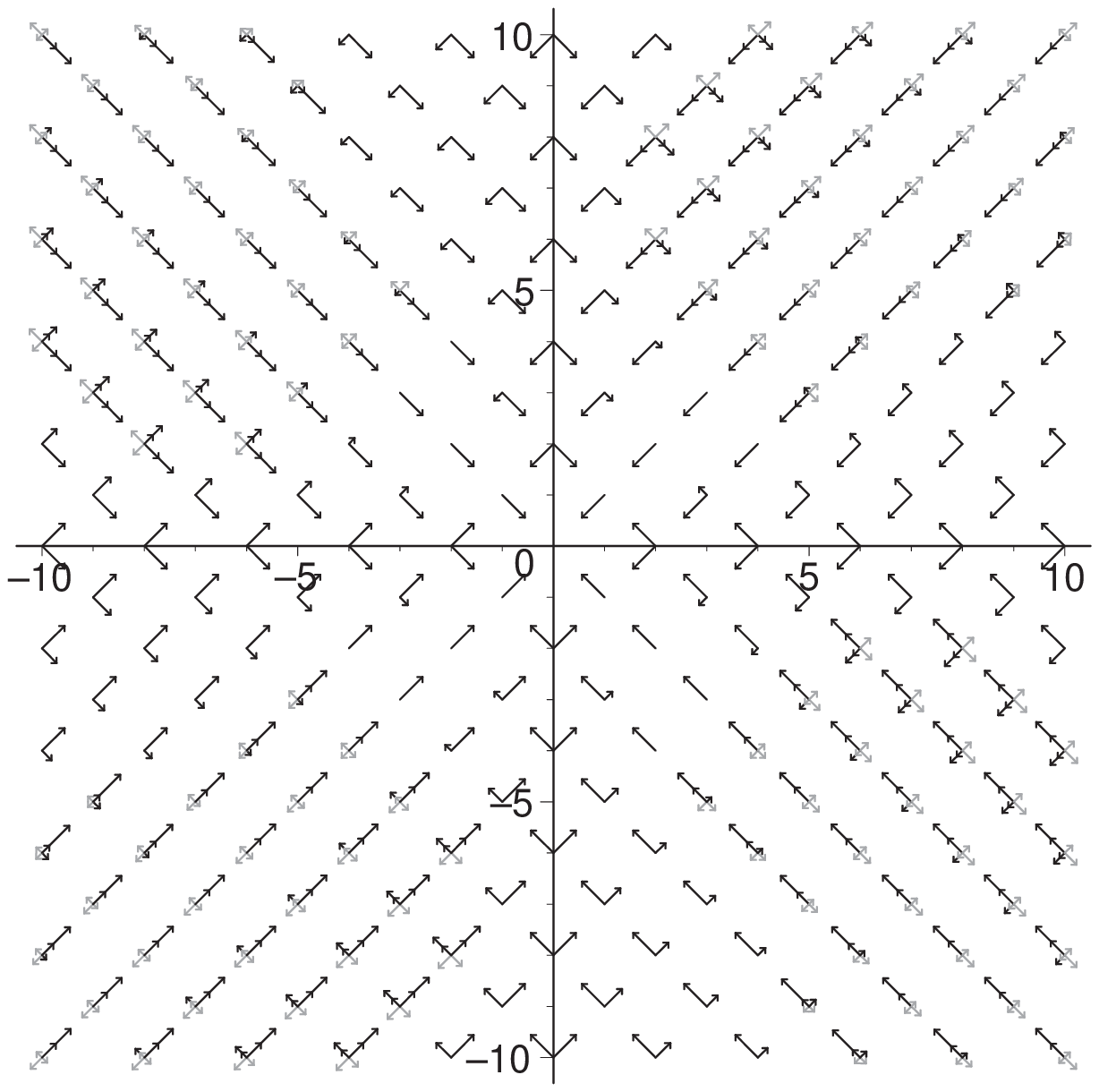}
        \label{subfig:plotdirs1}
    }
    \\
    \subfloat[Non-circular rotor sequence $(\arrpp,\arrmp,\arrpm,\arrmm)$.]{
        \includegraphics[bb=129pt 221pt 480pt 570pt,width=\ifCPC.65\else.6\fi\textwidth,clip]{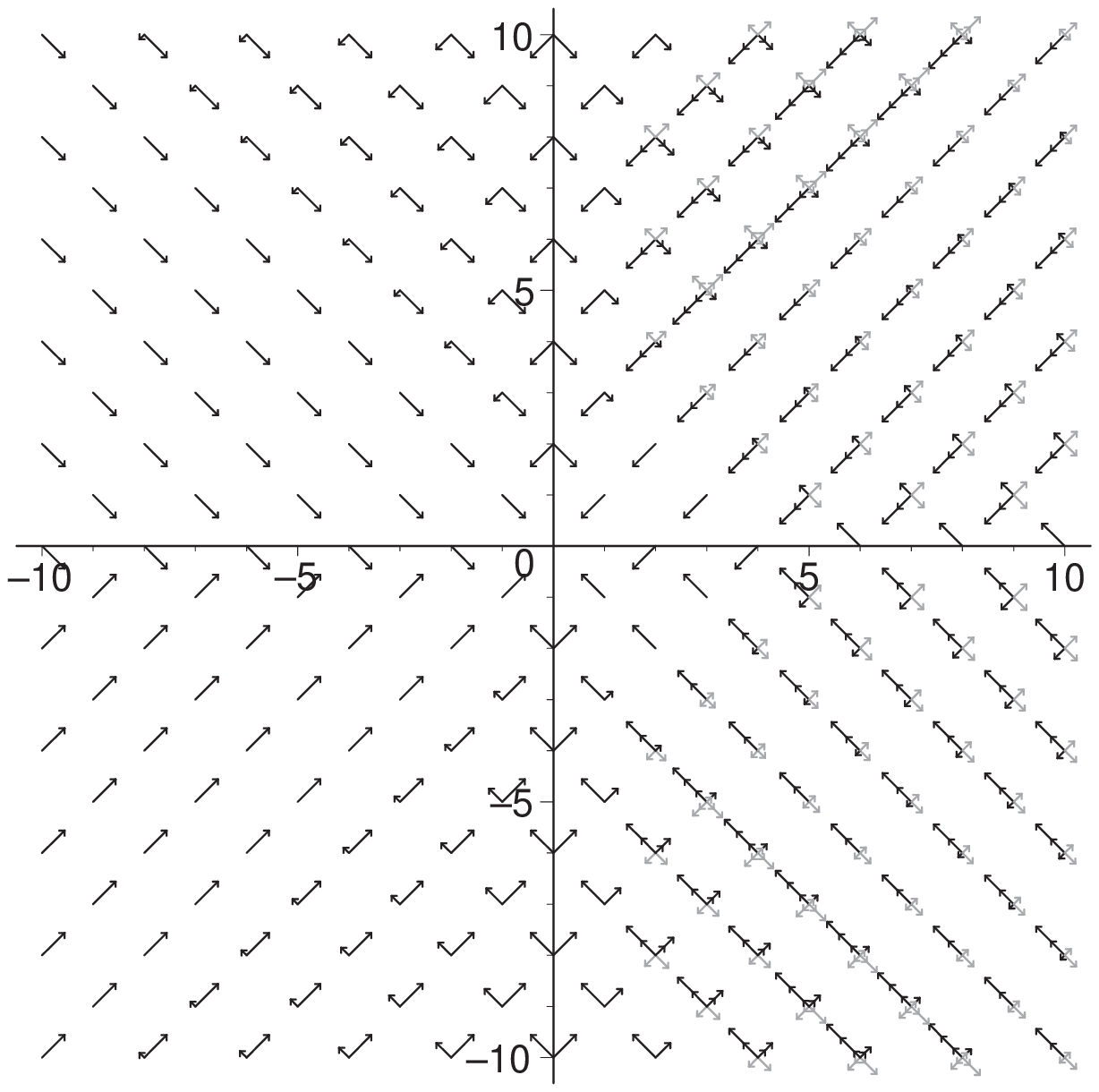}
        \label{subfig:plotdirs2}
    }
    \caption{At each position, the directions in which odd chips are sent in a
             block configuration are shown.
             If there are multiple chips sent from one position, the arrows are scaled
             according to their relative contribution.
             Grey arrows represent negative contributions.
             For position $\tbinom{5}{9}$, the second and third row of the second table on page~\pageref{tab:rotors}
             contains more details on the times when these chips are sent.}
    \label{fig:plotdirs}
\end{figure}

\begin{figure}[p]
    \centering
    \subfloat[Counterclockwise rotor sequence $(\arrpp,\arrmp,\arrmm,\arrpm)$.]{
        \includegraphics[bb=129pt 221pt 480pt 570pt,width=\ifCPC.65\else.6\fi\textwidth,clip]{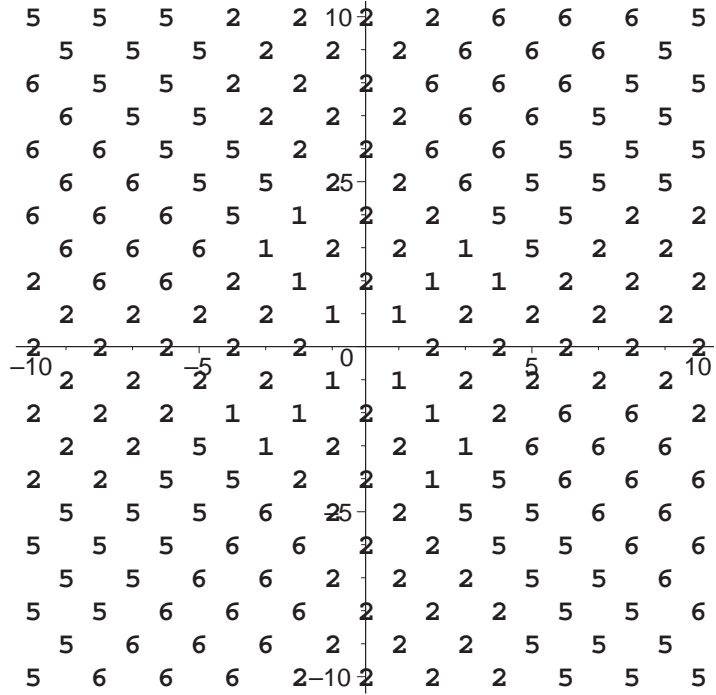}
        \label{subfig:numchips1}
    }
    \\
    \subfloat[Non-circular rotor sequence $(\arrpp,\arrmp,\arrpm,\arrmm)$.]{
        \includegraphics[bb=129pt 221pt 480pt 570pt,width=\ifCPC.65\else.6\fi\textwidth,clip]{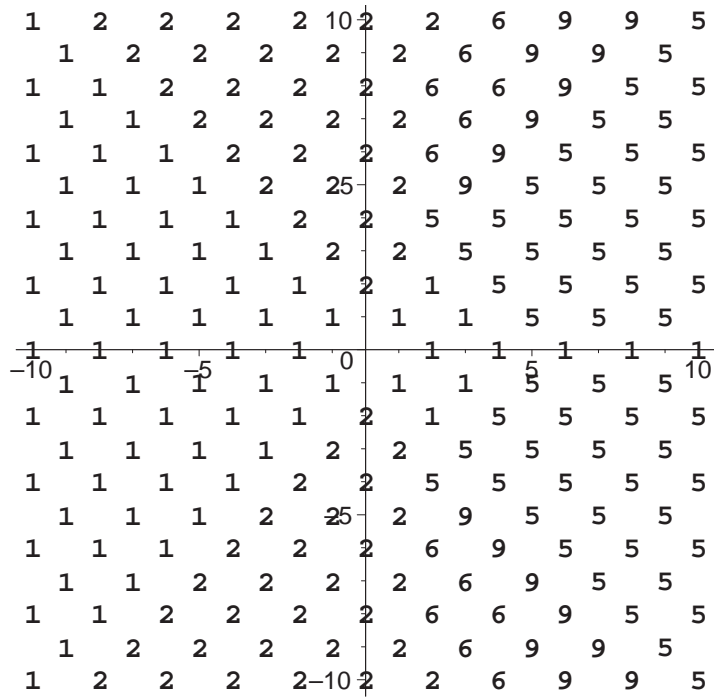}
        \label{subfig:numchips2}
    }
    \caption{Number of odd chips sent in a block configuration.  The 
             largest number of nine odd chips occurs only in \protect\subref{subfig:numchips2},
             not in \protect\subref{subfig:numchips1}.}
    \label{fig:numchips}
\end{figure}

\begin{figure}[p]
    \centering
    \subfloat[Counterclockwise rotor sequence $(\arrpp,\arrmp,\arrmm,\arrpm)$.]{
        \includegraphics[bb=129pt 221pt 480pt 570pt,width=\ifCPC.65\else.6\fi\textwidth
        								,clip]{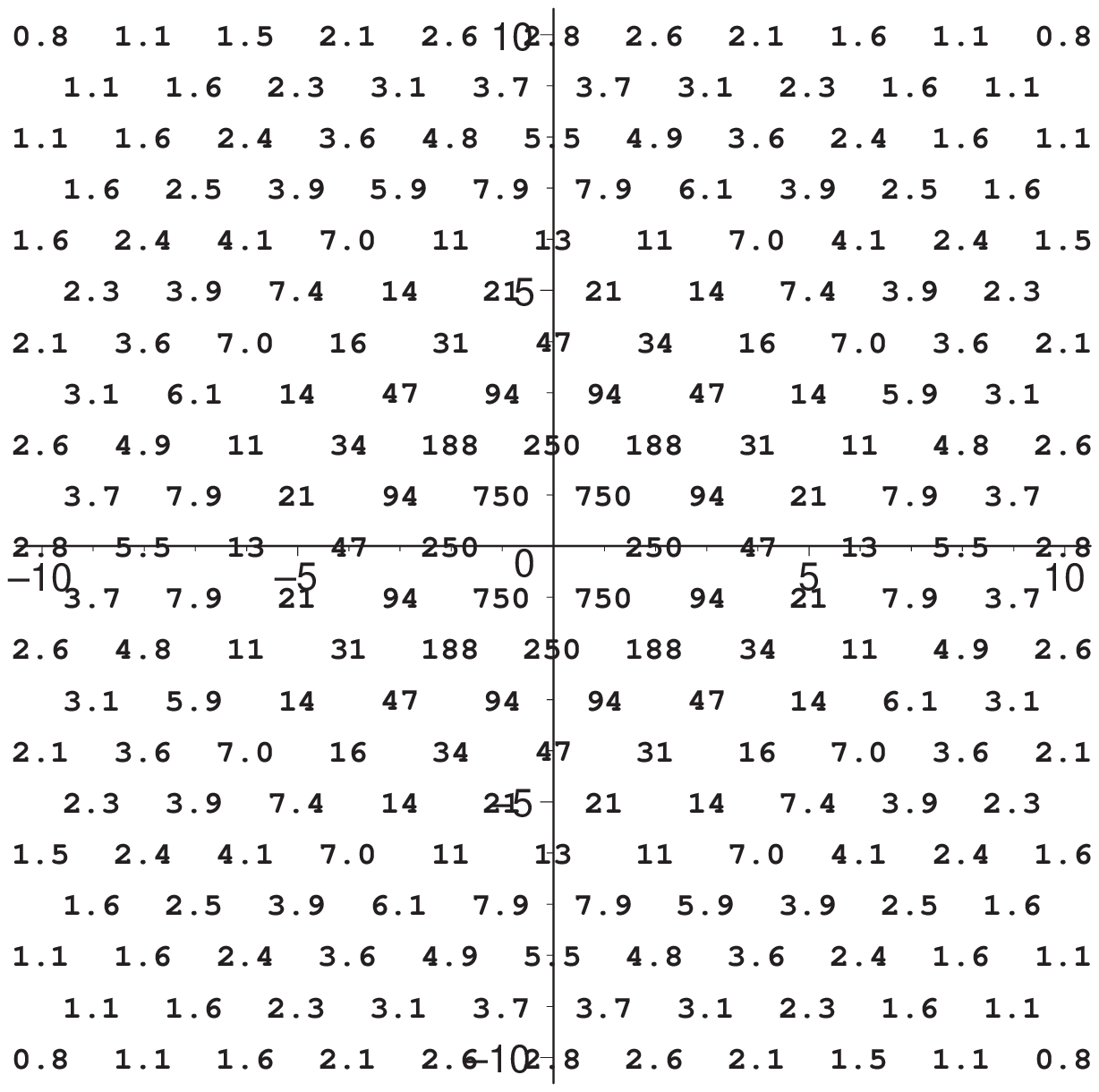}
        \label{subfig:val1}
    }
    \\
    \subfloat[Non-circular rotor sequence $(\arrpp,\arrmp,\arrpm,\arrmm)$.]{
        \includegraphics[bb=129pt 221pt 480pt 570pt,width=\ifCPC.65\else.6\fi\textwidth,clip]{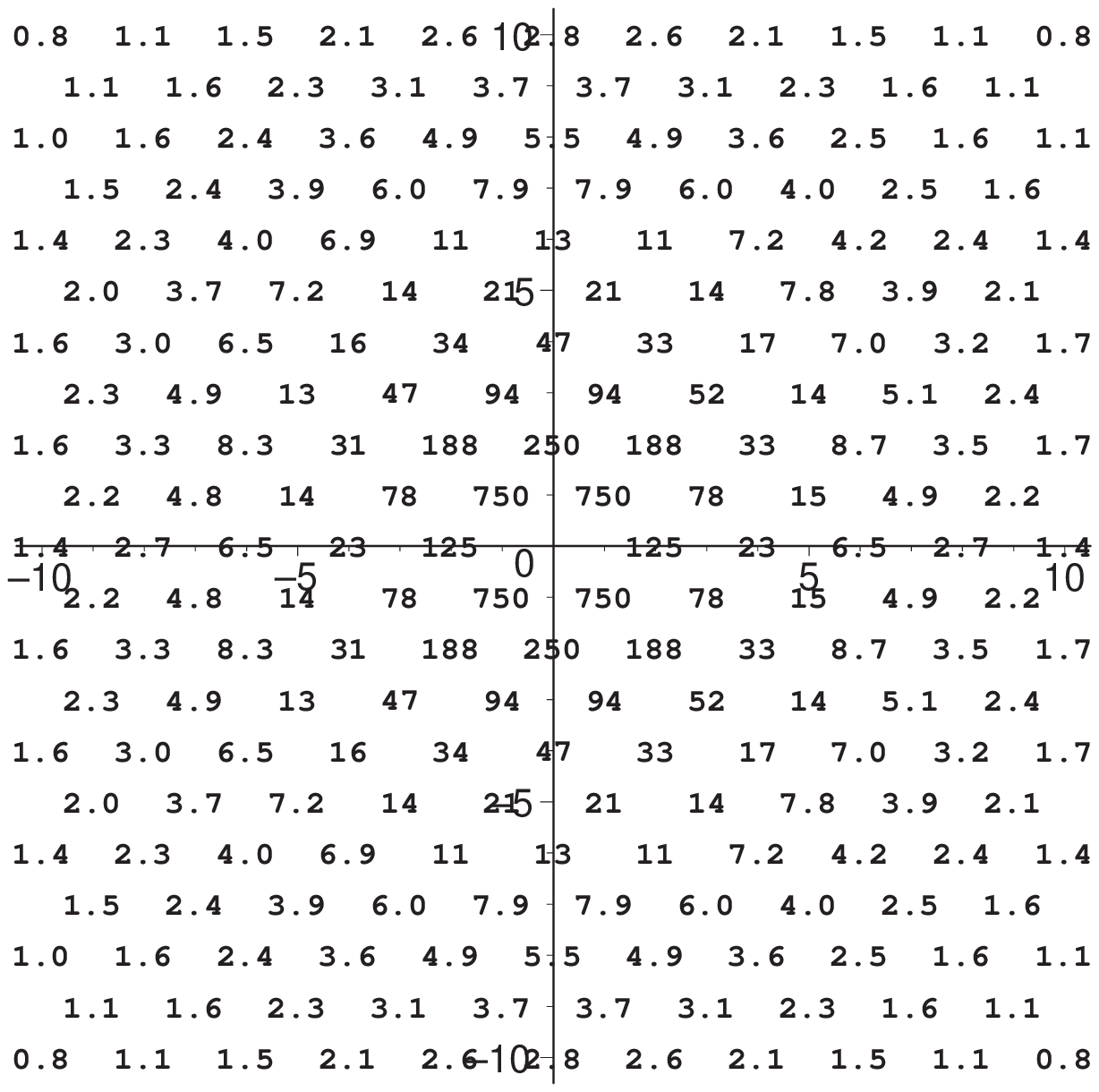}
        \label{subfig:val2}
    }
    \caption{Contributions $\MAXCON(\x)$ times 1000.
             Note the symmetries discussed in \secref{upper} and
             the quick descent of $\MAXCON(\x)$ for increasing $\|\x\|_2$.}
    \label{fig:val}
    \ifCPC\label{lastpage}\fi
\end{figure}

In Figures~\ref{fig:plotdirs} 
to~\ref{fig:val}, more information about the behavior of different positions in 
a worst-case setting (aiming at a surplus of chips at the origin) is collected. 

We also brief{}ly regarded the IDLA aggregation model. We saw that the 
surprisingly strong convergence to a Eucledian ball observed in earlier research 
also holds for non-circular rotor sequences and non-regular initial rotor 
settings. However, the suspected constant again seems to depend on the rotor 
sequences, and again, the circular ones seem to behave slightly worse than the non-circular ones.


\section*{Acknowledgments}

We would like to thank Joel Spencer and Jim Propp for several very inspiring discussions.


\ifCPC 
    \def\newblock{\hskip .11em plus .33em minus .07em}
\fi

\end{document}